\documentclass[12pt]{article}
\usepackage{latexsym, amssymb}
\usepackage{mathbbol}

\DeclareMathAlphabet\gothic{U}{euf}{m}{n}

\makeatletter
\def\eqnarray{\stepcounter{equation}\let\@currentlabel=\theequation
\global\@eqnswtrue
\tabskip\@centering\let\\=\@eqncr
$$\halign to \displaywidth\bgroup\hfil\global\@eqcnt\z@
  $\displaystyle\tabskip\z@{##}$&\global\@eqcnt\@ne
  \hfil$\displaystyle{{}##{}}$\hfil
  &\global\@eqcnt\tw@ $\displaystyle{##}$\hfil
  \tabskip\@centering&\llap{##}\tabskip\z@\cr}

\def\endeqnarray{\@@eqncr\egroup
      \global\advance\c@equation\m@ne$$\global\@ignoretrue}
\makeatother

\begin{document}
\bibliographystyle{tom}

\newtheorem{lemma}{Lemma}[section]
\newtheorem{thm}[lemma]{Theorem}
\newtheorem{cor}[lemma]{Corollary}
\newtheorem{voorb}[lemma]{Example}
\newtheorem{rem}[lemma]{Remark}
\newtheorem{prop}[lemma]{Proposition}
\newtheorem{stat}[lemma]{{\hspace{-5pt}}}

\newenvironment{remarkn}{\begin{rem} \rm}{\end{rem}}
\newenvironment{exam}{\begin{voorb} \rm}{\end{voorb}}

\newcounter{teller}
\renewcommand{\theteller}{\Roman{teller}}
\newenvironment{tabel}{\begin{list}%
{\rm \bf \Roman{teller}.\hfill}{\usecounter{teller} \leftmargin=1.1cm
\labelwidth=1.1cm \labelsep=0cm \parsep=0cm}
                      }{\end{list}}

\newcommand{\Ni}{{\bf N}}
\newcommand{\Ri}{{\bf R}}
\newcommand{\Ci}{{\bf C}}
\newcommand{\Ti}{{\bf T}}
\newcommand{\Zi}{{\bf Z}}
\newcommand{\Fi}{{\bf F}}

\newcommand{\proof}{\mbox{\bf Proof} \hspace{5pt}} 
\newcommand{\remark}{\mbox{\bf Remark} \hspace{5pt}}
\newcommand{\ruimte}{\vskip10.0pt plus 4.0pt minus 6.0pt}

\newcommand{\supp}{\mathop{\rm supp}}
\newcommand{\essinf}{\mathop{\rm ess\,inf}}
\newcommand{\esssup}{\mathop{\rm ess\,sup}}
\newcommand{\one}{\mathbb{1}}

\hyphenation{groups}
\hyphenation{unitary}

\newcommand{\cc}{{\cal C}}
\newcommand{\cd}{{\cal D}}
\newcommand{\ce}{{\cal E}}
\newcommand{\cf}{{\cal F}}
\newcommand{\ch}{{\cal H}}
\newcommand{\ci}{{\cal I}}
\newcommand{\ck}{{\cal K}}
\newcommand{\cl}{{\cal L}}
\newcommand{\cm}{{\cal M}}
\newcommand{\co}{{\cal O}}
\newcommand{\cs}{{\cal S}}
\newcommand{\ct}{{\cal T}}
\newcommand{\cx}{{\cal X}}
\newcommand{\cy}{{\cal Y}}
\newcommand{\cz}{{\cal Z}}

\thispagestyle{empty}

\begin{center}
{\Large\bf Dirichlet forms and degenerate elliptic operators  } \\[5mm]
\large A.F.M. ter Elst$^1$, Derek W. Robinson$^2$, Adam Sikora$^3$
and Yueping Zhu$^4$
\end{center}

\vspace{5mm}

\begin{center}
{\em Dedicated to Philippe Clement on the occasion of his retirement}
\end{center}

\vspace{5mm}

\begin{center}
{\bf Abstract}
\end{center}

\begin{list}{}{\leftmargin=1.8cm \rightmargin=1.8cm \listparindent=10mm 
   \parsep=0pt}
\item
It is shown that the theory of real symmetric second-order elliptic operators
in divergence form on $\Ri^d$ can be formulated in terms of a regular strongly local Dirichlet form
irregardless of the order of degeneracy.
The behaviour of the corresponding evolution semigroup $S_t$ can be described in terms of a
function $(A,B) \mapsto d(A\,;B)\in[0,\infty]$ over pairs of measurable subsets of $\Ri^d$. 
Then
\[
|(\varphi_A,S_t\varphi_B)|\leq e^{-d(A;B)^2(4t)^{-1}}\|\varphi_A\|_2\|\varphi_B\|_2
\]
for all $t>0$ and all $\varphi_A\in L_2(A)$, $\varphi_B\in L_2(B)$.
Moreover $S_tL_2(A)\subseteq L_2(A)$ for all $t>0$ if and only if $d(A\,;A^c)=\infty$ where $A^c$ denotes
the complement of $A$.
\end{list}

\vspace{1cm}
\noindent
September 2005 

\vspace{5mm}

\vspace{5mm}

\noindent
{\bf Home institutions:}    \\[3mm]
\begin{tabular}{@{}cl@{\hspace{10mm}}cl}
1. & Department of Mathematics  & 
  2. & Centre for Mathematics   \\
& \hspace{15mm} and Computing Science & 
  & \hspace{15mm} and its Applications  \\
& Eindhoven University of Technology & 
  & Mathematical Sciences Institute  \\
& P.O. Box 513 & 
  & Australian National University  \\
& 5600 MB Eindhoven & 
  & Canberra, ACT 0200  \\
& The Netherlands & 
  & Australia  \\[8mm]
3. & Department of Mathematical Sciences & 
 4. &  Department of Mathematics   \\
& New Mexico State University & 
  & Nantong University  \\
& P.O. Box 30001 & 
  & Nantong, 226007   \\
& Las Cruces & 
  & Jiangsu Province  \\
& NM 88003-8001, USA & 
  & P.R. China  \\
& {} & 
  & {} 
\end{tabular}

\newpage
\setcounter{page}{1}

\section{Introduction}\label{Sdfd1}

The usual starting point for the analysis of second-order divergence-form elliptic operators
with measurable coefficients is the precise definition of the operator by  quadratic form
techniques.
Let $c_{ij}=c_{ji}\in L_\infty(\Ri^d\!:\!\Ri)$, the  real-valued bounded measurable functions on $\Ri^d$,
and assume  that the 
$d\times d$-matrix $C=(c_{ij})$ is positive-definite almost-everywhere.
Then define the quadratic form $h$ by
\begin{equation}
h(\varphi)=\sum^d_{i,j=1}(\partial_i\varphi,c_{ij}\partial_j\varphi)
\label{edfd1.1}
\end{equation}
where $\partial_i=\partial/\partial x_i$ and $\varphi\in D(h)=W^{1,2}(\Ri^d)$.
It follows that $h$ is positive and the corresponding sesquilinear form $h(\cdot\,,\cdot) $ 
is symmetric.
Therefore if $h$ is closed there is a canonical construction which gives a unique positive
self-adjoint operator $H$ such that $D(h)=D(H^{1/2})$, $h(\varphi)=\|H^{1/2}\varphi\|_2^2$
for all $\varphi\in D(h)$ and $h(\psi,\varphi)=(\psi, H\varphi)$ if $\psi\in D(h)$
and $\varphi\in D(H)$.
Here and in the sequel $\|\cdot\|_p$ denotes the $L_p$-norm.
Formally one has
\[
H=-\sum^d_{i,j=1}\partial_i\, c_{ij}\,\partial_j
\;\;\;.
\]
The critical point in the form  approach is that $h$ must be closed.
But $h$ is closed if and only if there is a $\mu>0$ such that $C\geq \mu I$
(see, for example, \cite{ERZ1}, Proposition~2), i.e.,
if and only if the operator $H$ is strongly elliptic.
Therefore the construction of $H$ is not directly applicable to degenerate elliptic operators.
Nevertheless refinements of the theory of quadratic forms allow a precise definition
of the elliptic operator.
There are two distinct cases.

First, it is possible that $h$ is closable on $W^{1,2}(\Ri^d)$.
Then one can  repeat the previous
construction to obtain the self-adjoint elliptic operator associated with the closure
$\overline{h}$ of~$h$.
For example, let $\Omega$ be an open subset of $\Ri^d$ with $|\partial \Omega| = 0$.
Suppose $\supp C=\overline \Omega$ and $C(x) \geq\mu I>0$ uniformly for all  $x \in \Omega$.
Then $h$ is closable.
Indeed the restriction of $h$ to $W^{1,2}(\Omega)$ is closed as a form
on the subspace $L_2(\Omega)$  and the corresponding self-adjoint
operator $H_\Omega$ can be interpreted as the strongly elliptic operator with coefficients $C$
acting on $L_2(\Omega)$ with Neumann boundary conditions.
Then the  form corresponding to the operator  $H_\Omega\oplus 0$ on $L_2(\Ri^d)$,
where $H_\Omega$ acts on the subspace
$L_2(\Omega)$, is the closure of $h$.

Secondly, it is possible that $h$ is not closable.
There are many sufficient conditions for closability 
of the form 
(\ref{edfd1.1}) (see, for example, \cite{FOT}, Section~3.1,
or \cite{MR}, Chapter~II) and if
 $d=1$ then necessary and sufficient conditions  are given  by \cite{FOT}, Theorem~3.1.6.
The latter conditions restrict the possible degeneracy.
But Simon~\cite{bSim5}, Theorems~2.1 and 2.2, has shown that  a general positive 
quadratic form $h$ 
 can be decomposed as a sum $h=h_r+h_s$ of two positive forms with 
$D(h_r)=D(h)=D(h_s)$ with
$h_r$  the largest closable form majorized by $h$.
Simon refers to $h_r$ as the regular part of $h$.
Then one can construct the operator $H_0$ associated with the closure
$h_0=\overline{h_r}$ of the regular part of $h$.
Note that if $h$ is closed then $h_0=h$ and if $h$ is closable then $h_0=\overline h$
so in both cases $H_0$ coincides with the previously defined elliptic operator.
Therefore this method can be considered as the generic method of defining the operator
associated with the form (\ref{edfd1.1}).

There is an alternative method of constructing $H_0$ which demonstrates more clearly its
relationship with the formal definition of the elliptic operator.
Introduce the form $l$ by $D(l)=D(h)=W^{1,2}(\Ri^d)$ and 
\[
l(\varphi)=\sum^d_{i=1}\|\partial_i\varphi\|_2^2
\;\;\;.
\]
Then $l$ is closed and the corresponding positive self-adjoint operator is the usual
Laplacian $\Delta$.
Furthermore for each $\varepsilon>0$ the form $h_\varepsilon=h+\varepsilon \,l$ with domain
$D(h_\varepsilon)=D(h)$ is closed.
The corresponding positive self-adjoint operators $H_\varepsilon$ 
are strongly elliptic operators with coefficients $c_{ij}+\varepsilon\,\delta_{ij}$. 
These operators
form a decreasing sequence which, by a result of Kato \cite{Kat1}, Theorem~VIII.3.11,
  converges in the strong
resolvent sense to a positive self-adjoint operator.
The limit  is the operator $H_0$ associated with the form $h_0$ by \cite{bSim5},
Theorem~3.2.
The latter method of construction of $H_0$ as a limit of the $H_\varepsilon$  was used in
\cite{ERZ1} and \cite{ERSZ1} and motivated the terminology viscosity operator for $H_0$
and viscosity form for $h_0$.
The form $h_r$  constructed by Simon also occurs in the context of nonlinear phenomena and
discontinuous media and is variously described as the regularization or relaxation
of $h$ (see, for example,  \cite{Bra} \cite {ET} \cite{Jos} \cite{DalM}  \cite{Mosco}
and references therein).

Although the viscosity form provides a basis for the analysis of degenerate elliptic operators
its indirect definition makes it difficult to deduce even the most straightforward properties.
Simon, \cite{bSim4} Theorem~2, proved that $D(h_0)$ consists of those $\varphi\in L_2(\Ri^d)$
for which there is a sequence $\varphi_n\in D(h)$ such that $\lim_{n \to \infty} \varphi_n = \varphi$
in $L_2$ and $\liminf_{n\to\infty} h(\varphi_n)<\infty$.
Moreover, $h_0(\varphi)$ equals the minimum of all $\liminf_{n \to \infty} h(\varphi_n)$, 
where the minimum is taken over all $\varphi_1,\varphi_2,\ldots \in D(h)$ such that 
$\lim_{n \to \infty} \varphi_n = \varphi$ in $L_2$.
(See \cite{bSim4}, Theorem~3.)
It was observed in \cite{ERZ1} that this characterization allows one to deduce that $h_0$ is
a Dirichlet form.
The first purpose of this note is to use the theory of Dirichlet forms (see \cite{FOT} 
\cite{BH} for background material) to strengthen the earlier conclusion.
Specifically we establish the following result in Section~\ref{Sdfd2}.

\begin{thm}\label{tdfd1.0}
The viscosity form $h_0$ is a
regular local Dirichlet form.
\end{thm}

The regularity is straightforward but the locality is not so evident and requires some control
of the dissipativity of  the semigroup $S^{(0)}$ generated by the viscosity operator $H_0$.
Note that we adopt the terminology of \cite{BH} which differs from that of \cite{FOT}.
Locality in \cite{BH} corresponds to strong locality in \cite{FOT} if the form is regular.
(See \cite{Schm} for a detailed study of the different forms of locality (in the sense of \cite{FOT}).)

Our second purpose is to extend earlier results on $L_2$ off-diagonal bounds for the semigroup
$S^{(0)}$ generated by the viscosity operator $H_0$.
These bounds, which are also referred to as Davies--Gaffney estimates, integrated Gaussian 
estimates or an integrated maximum principle (see, for example, 
\cite{Aus2} \cite{CGT}  \cite{Dav12} \cite{Gri3}  \cite{Stu4} \cite{Stu2}),
give upper bounds on the cross-norm $\|S^{(0)}_t\|_{A\to B}$ of the semigroup  $S^{(0)}_t$ between
 subspaces $L_2(A)$ and $L_2(B)$ of the form
\begin{equation}
\|S^{(0)}_t\|_{A\to B}\leq e^{-d(A;B)^2(4t)^{-1}}
\label{edfd237}
\end{equation}
where $d(A\,;B)$ is an appropriate  measure of distance between the subsets $A$ and $B$ of $\Ri^d$.
In the case of strongly elliptic operators there is an essentially unique distance,
variously referred to as the Riemannian distance, the intrinsic distance or the control
distance, suited to the problem.
But for degenerate elliptic operators, which exhibit phenomena of separation and isolation
\cite{ERSZ1}, the Riemannian distance is not necessarily appropriate.

We  will establish  estimates in terms of a set-theoretic function which is not strictly a 
distance since it can take the value infinity.
The function is defined by a version of a standard
variational principle which has been used widely in the analysis of strongly elliptic 
operators and which was extended to the general theory of Dirichlet forms by Biroli and Mosco
\cite{BM1} \cite{BM}.
In the case of strongly elliptic operators or subelliptic operators with smooth coefficients
this method gives a distance equivalent to the Riemannian distance obtained by path methods
(see, for example, \cite{JSC}, Section~3).
In the degenerate situation we make a choice of the set of variational functions which gives
a distance compatible with the separation properties.
In particular $d(A\,;A^c)=\infty$ if and only if  $S^{(0)}_tL_2(A)\subseteq L_2(A)$,
where $A^c$ is the complement of $A$.
Our choice is aimed to maximize the distance and thereby optimize the bounds~(\ref{edfd237}).
We carry out the analysis in a  general setting of local Dirichlet forms
which is applicable to degenerate elliptic operators defined as above but has a wider 
range of applicability.

Let $X$ be a topological Hausdorff space equipped with a $\sigma$-finite Borel measure $\mu$.
Further let $\ce$ be a local Dirichlet form on $L_2(X)$ in the sense of \cite{BH}.
First, 
for all $\psi \in D(\ce) \cap L_\infty(X\!:\!\Ri)$ define 
$\ci^{(\ce)}_\psi \colon D(\ce) \cap L_\infty(X\!:\!\Ri) \to \Ri$ by
\[
\ci^{(\ce)}_\psi(\varphi) = \ce(\psi \, \varphi,\psi) - 2^{-1} \ce(\psi^2,\varphi)
\;\;\; .  \]
If no confusion is possible  we drop the suffix and write 
$\ci_\psi(\varphi) = \ci^{(\ce)}_\psi(\varphi)$.
If $\varphi\geq 0$ it follows that 
$\psi \mapsto \ci_\psi(\varphi)\in \Ri$
is a Markovian form with  domain $D(\ce)\cap L_\infty(X\!:\!\Ri)$
(see  \cite{BH}, Proposition~I.4.1.1).
This form  is referred to as the truncated form by Roth
 \cite{Roth}, Theorem~5.

Secondly, define $D(\ce)_{\rm loc}$ as the vector space of 
(equivalent classes of) all measurable functions $\psi \colon X \to \Ci$
such that for every compact subset $K$ of $X$ there is a $\hat \psi \in D(\ce)$
with $\psi|_K = \hat \psi|_K$.
Since $\ce$ is  local one can define 
$\widehat \ci_\psi^{(\ce)} = \widehat \ci_\psi \colon D(\ce)\cap L_{\infty,c}(X\!:\!\Ri) \to \Ri$ by
\[
\widehat \ci_\psi(\varphi) = \ci_{\hat \psi}(\varphi)
\]
for all $\psi \in D(\ce)_{\rm loc} \cap L_\infty(X\!:\!\Ri)$ and 
$\varphi \in D(\ce) \cap L_{\infty,c}(X\!:\!\Ri)$ where 
$\hat \psi \in D(\ce)\cap L_\infty(X\!:\!\Ri)$ is such that $\psi|_{\supp \varphi} = \hat \psi|_{\supp \varphi}$.
Here $L_{\infty,c}(X\!:\!\Ri) = \{ \varphi \in L_\infty(X\!:\!\Ri) : \supp \varphi \mbox{ is compact} \} $
and $(\supp \varphi)^{\rm c}$ is the union of all open subsets $U \subseteq X$ such that
$\varphi|_U = 0$ almost everywhere.
Thirdly, for all $\psi \in D(\ce)_{\rm loc} \cap L_\infty(X\!:\!\Ri)$ define
\[
|||\widehat \ci_\psi|||
= \sup \{ |\widehat \ci_\psi(\varphi)| : 
      \varphi \in D(\ce) \cap L_{\infty,c}(X\!:\!\Ri) , \; \|\varphi\|_1 \leq 1 \}
\in [0,\infty]
\;\;\; .  \]
Fourthly, for all $\psi \in L_\infty(X\!:\!\Ri)$ and measurable sets $A,B \subset X$ introduce 
\begin{eqnarray*}
d_\psi(A\,;B)
&=& \sup \{ M \in \Ri : \psi(a) - \psi(b) \geq M \mbox{ for a.e.\ $a \in A$ and a.e.\ }
                      b \in B \}  \\[5pt]
&=&\essinf_{x\in A}\psi(x)-\esssup_{y\in B}\psi(y) 
\in \langle -\infty,\infty]
\;\;\; .  
\end{eqnarray*}
(Recall that $\esssup_{y\in B}\psi(y) = \inf \{ m \in \Ri : | \{ y \in B : \psi(y) > m \} | = 0 \} \in [-\infty,\infty \rangle$
and $\essinf_{x \in A} \psi(x) = - \esssup_{x \in A} -\psi(x)$.) 
Finally define the set theoretic distance
\[
d(A\,;B)
= d^{(\ce)}(A\,;B)
= \sup \{ d_\psi(A\,;B) : \psi \in D_0(\ce)  \}
\;\;\;,
\]
where
\[
D_0(\ce)
= \{ \psi\in D(\ce)_{\rm loc} \cap L_\infty(X\!:\!\Ri): |||\widehat \ci_\psi||| \leq 1 \}
\;\;\;.
\]
A similar definition was given by Hino and Ramirez \cite{HiR},
but since they consider probability spaces the introduction of $D(\ce)_{\rm loc}$ is 
unnecessary in \cite{HiR}.
If, however, we were to replace $D(\ce)_{\rm loc}$ by $D(\ce)$ in the definition 
of $d(A\,;B)$ then, since $\psi\in L_2(X)$, one would obtain $d_\psi(A\,;B)\leq  0$ for all measurable sets $A,B \subset X$
with $|A| = |B| = \infty$ and this would give $d(A\,;B)= 0$.
On the other hand  the definition with $D(\ce)_{\rm loc}$ is not useful unless
$D(\ce)$ contains sufficient bounded functions of compact support.

\begin{thm} \label{tdfd3.1}
Let $\ce$ be a    local Dirichlet form  on $L_2(X)$ with $\one \in D(\ce)_{\rm loc}$ and 
such that $D(\ce) \cap L_{\infty,c}(X)$ is a core for $\ce$.
Further let $d(\cdot\,;\cdot)$ denote the corresponding
set-theoretic distance.
If  $S$ denotes the semigroup generated by the self-adjoint  operator $H$
on $L_2(X)$ associated with $\ce$ and if 
 $A$ and $B$ are measurable subsets of $X$ then
\[
|(\varphi_A, S_t\varphi_B)|\leq e^{-d(A\,;B)^2(4t)^{-1}}\|\varphi_A\|_2\|\varphi_B\|_2
\]
for all $\varphi_A\in L_2(A)$, $\varphi_B\in L_2(B)$ and $t>0$ with the convention
$e^{-\infty}=0$.
\end{thm}

This applies in particular to the viscosity operator $H_0$.
One has $D(l) = D(h) \subseteq D(h_0)$.
Hence $\one \in D(l)_{\rm loc} \subseteq D(h_0)_{\rm loc}$.

Theorem~\ref{tdfd3.1} will be proved in Section~\ref{Sdfd3}.

If $h$ is closed, i.e., if the corresponding operator is strongly elliptic, then $d^{(h)}(\cdot\,;\cdot)$
is finite-valued.
Nevertheless for degenerate operators one can have $d^{(h)}(A\,;B)=\infty$ with $A,B$ non-empty
open subsets of $\Ri^d$ and the action of the semigroup 
$S^{(0)}$ can be non-ergodic \cite{ERSZ1}.

Specifically we establish the following result in Section~\ref{Sdfd4}.

\begin{thm} \label{tdfd401}
Let $\ce$ be a    local Dirichlet form  on $L_2(X)$ with $\one \in D(\ce)_{\rm loc}$ and 
such that $D(\ce) \cap L_{\infty,c}(X)$ is a core for $\ce$.
Further let $d(\cdot\,;\cdot)$ denote the corresponding
set-theoretic distance.
If  $S$ denotes the semigroup generated by the self-adjoint  operator $H$
on $L_2(X)$ associated with $\ce$
and $A\subset X$ is measurable then the following conditions are 
equivalent.
\begin{tabel}
\item \label{tdfd401-0}
$S_t L_2(A) \subseteq L_2(A)$ for one $t > 0$.
\item \label{tdfd401-1}
$S_t L_2(A) \subseteq L_2(A)$ for all $t > 0$.
\item \label{tdfd401-2}
$d(A\,;A^{\rm c}) = \infty$.
\item \label{tdfd401-4}
$d(A\,;A^{\rm c}) > 0$.
\end{tabel}
\end{thm}

We conclude by deriving alternate characterizations of the distance and deriving some of 
its general properties in Section~\ref{Sdfd5}.

\section{Dirichlet forms} \label{Sdfd2}

In this section we prove that the form $h_0$ defined in the introduction 
is a regular strongly local Dirichlet form.
First we recall the basic definitions.

A Dirichlet form $\ce$ on $L_2(X)$ is called {\bf regular} if there is a subset of $D(\ce) \cap C_c(X)$
which is  a core of $\ce$, i.e., which is dense  in $D(\ce)$
with respect to the natural norm $\varphi \mapsto (\ce(\varphi) + \|\varphi\|_2^2)^{1/2}$,
and which is also dense  in $C_0(X)$ with respect to the supremum norm.
There are three kinds of locality for Dirichlet forms \cite{BH} \cite{FOT}.
For locality we choose the definition of \cite{BH}.
A Dirichlet form $\ce$ is called {\bf local} if $\ce(\psi,\varphi) = 0$ for all $\varphi,\psi \in D(\ce)$
and $a \in \Ri$ such that $(\varphi + a \one) \psi = 0$.
Alternatively, the Dirichlet form $\ce$ is called {\bf \cite{FOT}-local} if $\ce(\psi,\varphi) = 0$ for all $\varphi,\psi \in D(\ce)$
with $\supp \varphi$ and $\supp \psi$ compact and $\supp \varphi \cap \supp \psi = \emptyset$.
Moreover, it is called {\bf \cite{FOT}-strongly local}
if $\ce(\psi,\varphi) = 0$ for all $\varphi,\psi \in D(\ce)$ with $\supp \varphi$ and 
$\supp \psi$ compact and $\psi$  constant on a neighbourhood of $\supp \varphi$.
Every \cite{FOT}-strongly local Dirichlet form is \cite{FOT}-local and if $X$ satisfies the
second axiom of countability
then every local Dirichlet form (in the sense of \cite{BH}) is \cite{FOT}-strongly local.
If $X$ is a locally compact separable
metric space and $\mu$ is a Radon measure such that  $\supp \mu = X$,
then a regular \cite{FOT}-strongly local Dirichlet form is local (in the sense of \cite{BH}) by \cite{BH}
Remark~I.5.1.5 and Proposition~I.5.1.3 $(L_0) \Rightarrow (L_2)$.

\begin{lemma} \label{ldirlocal3}
The form $h_0$ is regular.
Moreover, every core for the form $l$ of the Laplacian is a core for $h_0$.
\end{lemma}
\proof\
Let $h_r$ be the regular part of $h$ as in \cite{bSim5}.
Then $D(h_r) = D(h)$ and $h_0$ is the closure of $h_r$ by definition.
So $D(l) = D(h) = D(h_r)$ is a core for $h_0$.
Since $h_0\leq h \leq \|C\|\, l$,
 where  $\|C\|$ denotes the essential supremum of the matrix norms $\|C(x)\|$,
it follows that 
any core for $l$ is also a core for $h_0$.

Finally, since $C_c^\infty(\Ri^d) \subset W^{1,2}(\Ri^d) = D(h) \subset D(h_0)$
the set $C_c^\infty(\Ri^d) \subset D(h_0) \cap C_c(\Ri^d)$ is a core 
of $h_0$ and is dense in $C_0(\Ri^d)$.\hfill$\Box$

\ruimte

Next we examine the locality properties.

\begin{prop} \label{pdirlocal4}
The form $h_0$ is local.
\end{prop}
\proof\
First note that the \cite{FOT}-locality of $h_0$ is an easy consequence of Lemma~3.5 in \cite{ERSZ1}.
Fix  $\varphi,\psi\in D(h_0)$ with $\supp \varphi \cap \supp \psi = \emptyset$.
If $D$ is the Euclidean distance between the support of $\varphi$ and the support of $\psi$
then
\begin{eqnarray*}
|h_0(\psi,\varphi)|=\lim_{t\to0}t^{-1}|(\psi,S^{(0)}_t\varphi)|
\leq\lim_{t\to0}t^{-1} e^{-D^2 (4\|C\|t)^{-1} }\|\psi\|_2 \|\varphi\|_2=0
\end{eqnarray*}
where we have used $(\psi,\varphi)=0$ and the statement  of \cite{ERSZ1}, Lemma~3.5.
The proof of \cite{FOT}-strong locality is similar but depends on the estimates derived  
in the proof of Lemma~3.5.

Fix  $\varphi,\psi \in D(h_0)$ with $\supp \varphi$ and 
$\supp \psi$ compact and $\psi=1$   on a neighbourhood $U$ of $\supp \varphi$.
It suffices, by the remark preceding Lemma~\ref{ldirlocal3}, to prove that $h_0(\psi,\varphi)=0$. 
We may assume the support of $\psi$ is contained in the Euclidean ball $B_R$ centred at the origin
and with radius $R>0$.
Set $\chi_R=\one_{B_R}$.
Since $S^{(0)}_t\one=\one$, by Proposition~3.6 of \cite{ERSZ1},   one has
$(\psi,\varphi) = (\one, \varphi) = (\one,S^{(0)}_t \varphi)$.
Therefore 
\begin{equation}
t^{-1}(\psi,(I-S^{(0)}_t)\varphi)
= t^{-1}((\chi_R-\psi),S^{(0)}_t\varphi)
   - t^{-1}((\chi_R-\one),S^{(0)}_t\varphi)
\;\;\;.  
\label{edfd136}
\end{equation}
Now let $D$ denote the Euclidean distance from $\supp\varphi$ to $U^{\rm c}$.
It follows by assumption that $D>0$.
Then by  Lemma~3.5 of \cite{ERSZ1} there is a $c>0$ such that 
\begin{equation}
t^{-1}|((\chi_R-\psi),S^{(0)}_t\varphi)|
\leq t^{-1}e^{-D^2 (4\|C\|t)^{-1} }\|\chi_R-\psi\|_2\|\varphi\|_2
\leq c\,R^{d/2}\,t^{-1}e^{-D^2 (4\|C\|t)^{-1}}
\label{edfd137}
\end{equation}
uniformly for all large $R$ and all $t>0$.
The factor $R^{d/2}$ comes from the $L_2$-norm of $\chi_R$. 
Alternatively the estimate used in the proof of 
 Proposition~3.6 in \cite{ERSZ1} establishes that there are 
$a,b > 0$ such that 
\[
|((\one - \chi_R), S^{(0)}_t \varphi)|
\leq a \sum_{n=2}^\infty n^{d/2} R^{d/2} e^{-b n^2 R^2 t^{-1}} \|\varphi\|_2
\]
uniformly for all $R,t > 0$ such that 
$\supp \varphi \subset B_{2^{-1} R}$.
Next there is a $c' > 0$ such that 
$\sum_{n=2}^\infty n^{d/2} \, e^{-\alpha n^2} \leq c' \, \alpha^{-(d+2)/4}e^{-\alpha}$
uniformly for all $\alpha > 0$.
Hence there is a $c''>0$ such that 
\[
t^{-1}|((\one - \chi_R), S^{(0)}_t \varphi)|
\leq c'' \, R^{-1}t^{(d-2)/4}e^{-bR^2t^{-1}}
\;\;\;.
\]
Combining this with (\ref{edfd136}) and (\ref{edfd137}) one has 
\[
t^{-1}|(\psi,(I-S^{(0)}_t)\varphi)|
\leq c\,R^{d/2}t^{-1}e^{-D^2 (4\|C\|t)^{-1} } + c'' \, R^{-1}t^{(d-2)/4}e^{-bR^2t^{-1}}
\]
for all large $R$ and all $t>0$.
Taking the limit $t\to0$ 
establishes that 
$h_0(\psi,\varphi)=0$.
Since $h_0$ is regular it follows that $h_0$ is local.\hfill$\Box$

\section{$L_2$ off-diagonal bounds} \label{Sdfd3}

In this section we prove the $L_2$ off-diagonal bounds of Theorem~\ref{tdfd3.1}.
Initially we assume that $\ce$ is a  local Dirichlet form on $L_2(X)$ without assuming 
any kind of regularity.
The proof of the theorem follows by
standard reasoning based on an exponential perturbation technique
used by Gaffney~\cite{Gaf} and subsequently developed by Davies~\cite{Dav7} in the proof of pointwise Gaussian bounds.
It is essential to establish that the perturbation of the Dirichlet form is quadratic.
But this is a general consequence of  locality.
To exploit the latter property we use a result of  Andersson~\cite{And}
 and \cite{Roth}.

\begin{prop}\label{pdirlocal6}
Let $\ce$ be a  local Dirichlet form on $L_2(X)$
and $\varphi_1,\ldots,\varphi_n \in D(\ce) \cap L_\infty(X\!:\!\Ri)$.
Then for all $i,j \in \{ 1,\ldots,n \} $ there exists a unique real Radon measure 
$\sigma_{ij}^{(\ce,\varphi_1,\ldots,\varphi_n)} = \sigma_{ij}^{(\varphi_1,\ldots,\varphi_n)}$
on $\Ri^n$ 
such that $\sigma_{ij}^{(\varphi_1,\ldots,\varphi_n)} = \sigma_{ji}^{(\varphi_1,\ldots,\varphi_n)}$
for all $i,j \in \{ 1,\ldots,n \} $ and 
\begin{equation}
\ce(F_0(\varphi_1,\ldots,\varphi_n), G_0(\varphi_1,\ldots,\varphi_n))
= \sum_{i,j=1}^n \int_{\Ri^n} d\sigma_{ij}^{(\varphi_1,\ldots,\varphi_n)} \, 
      \overline{ \frac{\partial F}{\partial x_i} } \, \frac{\partial G}{\partial x_i}
\label{epdirlocal6;1}
\end{equation}
for all $F,G \in C^1(\Ri^n)$ where $F_0=F-F(0)$ and $G_0=G-G(0)$.
Let $K$ be a compact subset of $\Ri^n$ such that 
$(\varphi_1(x),\ldots,\varphi_n(x)) \in K$ for a.e.\ $x \in X$.
Then $\supp \sigma_{ij}^{(\varphi_1,\ldots,\varphi_n)} \subseteq K$ for all 
$i,j \in \{ 1,\ldots,n \} $.
In particular, if $i \in \{ 1,\ldots,n \} $ then $\sigma_{ii}^{(\varphi_1,\ldots,\varphi_n)}$ is a 
finite $($positive$)$ measure.

Moreover, if  $\cf$ is a second local Dirichlet form with $\ce \leq \cf$ then 
\begin{equation}
\int d\sigma_{ii}^{(\ce,\varphi_1,\ldots,\varphi_n)} \, \chi
\leq \int d\sigma_{ii}^{(\cf,\varphi_1,\ldots,\varphi_n)} \, \chi
\label{epdirlocal6;2}
\end{equation}
for all  $\varphi_1,\ldots,\varphi_n \in D(\cf) \cap L_\infty(X\!:\!\Ri)$, 
$i \in \{ 1,\ldots,n \} $ and $\chi \in C_c(\Ri^d)$ with $\chi \geq 0$.
\end{prop}
\proof\
 Theorem~I.5.2.1 of \cite{BH}, which elaborates a result of Andersson~\cite{And}, establishes that
there exist unique real Radon measures 
$\sigma_{ij}^{(\varphi_1,\ldots,\varphi_n)}$ such that 
$\sigma_{ij}^{(\varphi_1,\ldots,\varphi_n)} = \sigma_{ji}^{(\varphi_1,\ldots,\varphi_n)}$
for all $i,j \in \{ 1,\ldots,n \} $ and (\ref{epdirlocal6;1}) is valid for all $F,G \in C_c^1(\Ri^d)$.
Moreover, $\sum_{i,j=1}^n \xi_i \, \xi_j \, \sigma_{ij}^{(\varphi_1,\ldots,\varphi_n)}$
is a (positive) measure for all $\xi \in \Ri^n$.
Let $\xi \in \Ri^n$ and $\chi \in C_c^1(\Ri^n)$.
For all $\lambda > 0$ define $F_\lambda \in C_c^1(\Ri^n)$ by 
$F_\lambda(x) = e^{i \lambda x \cdot \xi} \, \chi(x) $.
Then it follows from (\ref{epdirlocal6;1}) that 
\begin{equation}
\lim_{\lambda \to \infty} 
    \lambda^{-2} \, \ce(F_{\lambda,0}(\varphi_1,\ldots,\varphi_n))
= \int \sum_{i,j=1}^n \xi_i \, \xi_j \, d\sigma_{ij}^{(\ce,\varphi_1,\ldots,\varphi_n)} \, |\chi|^2
\;\;\; .  
\label{epdirlocal6;3}
\end{equation}
So if $\supp \chi \subset K^{\rm c}$ then the left hand side of (\ref{epdirlocal6;3})
vanishes.
Hence $\supp \sigma_{ij}^{(\varphi_1,\ldots,\varphi_n)} \subseteq K$.
But then (\ref{epdirlocal6;1}) extends to all $F,G \in C^1(\Ri^n)$.
Moreover the measure $\sigma_{ii}^{(\varphi_1,\ldots,\varphi_n)}$ is finite for 
all $i \in \{ 1,\ldots,n \} $ since $\sigma_{ii}^{(\varphi_1,\ldots,\varphi_n)}$ is regular.

Finally, if  $\cf$ is a second  local Dirichlet form with $\ce \leq \cf$ then
it follows from (\ref{epdirlocal6;3}), applied to $\ce$ and $\cf$, that 
\[
\int \sum_{i,j=1}^n \xi_i \, \xi_j \, d\sigma_{ij}^{(\ce,\varphi_1,\ldots,\varphi_n)} \, |\chi|^2
\leq \int \sum_{i,j=1}^n \xi_i \, \xi_j \, d\sigma_{ij}^{(\cf,\varphi_1,\ldots,\varphi_n)} \, |\chi|^2
\]
for all $\chi \in C_c^1(\Ri^n)$ and $\xi \in \Ri^n$.
Setting $\xi = e_i$ gives (\ref{epdirlocal6;2}) by a density argument.\hfill$\Box$

\ruimte

Now we are prepared to prove the essential perturbation result for elements in $D(\ce)\cap L_\infty(X\!:\!\Ri)$.

\begin{prop}\label{pdirlocal5}
If $\ce$ is a  local Dirichlet form then
\[
\ce(\varphi,\varphi) - \ce(e^{-\psi}\varphi,e^{\psi}\varphi)
= \ci^{(\ce)}_\psi(\varphi^2)
\]
for all $\varphi,\psi\in D(\ce)\cap L_\infty(X\!:\!\Ri)$.
Moreover, if  $\cf$ is  a second  local Dirichlet form with $\ce \leq \cf$ then 
\[
\ci^{(\ce)}_\psi(\varphi)
\leq \ci^{(\cf)}_\psi(\varphi)
\]
for all $\varphi,\psi\in D(\cf) \cap L_\infty(X\!:\!\Ri)$ with $\varphi \geq 0$.
\end{prop}
\proof\
Observe that
\begin{eqnarray*}
\ce(\varphi,\varphi)-\ce(e^{-\psi}\varphi,e^{\psi}\varphi)
&=&-\ce((e^{-\psi}-1)\varphi,(e^{\psi}-1)\varphi)\\[5pt]
&&\hspace{2cm}{}-\ce((e^{-\psi}-1)\varphi,\varphi)
-\ce(\varphi,(e^{\psi}-1)\varphi)
\;\;\;.
\end{eqnarray*}
Now we can apply Proposition~\ref{pdirlocal6} with $n=2$ to each term on the right hand side.

First, one has
\begin{eqnarray*}
-\ce((e^{-\psi}-1)\varphi,(e^{\psi}-1)\varphi)
&=&\int d\sigma^{(\psi,\varphi)}_{1,1}(x_1,x_2)\,x_2^2\\[5pt]
&&\hspace{1cm}{}+\int d\sigma^{(\psi,\varphi)}_{1,2}(x_1,x_2)\,
(e^{-x_1}(e^{x_1}-1)-e^{x_1}(e^{-x_1}-1))x_2\\[5pt]
&&\hspace{2cm}{}-\int d\sigma^{(\psi,\varphi)}_{2,2}(x_1,x_2)\,(e^{-x_1}-1)(e^{x_1}-1)
\;\;\;.
\end{eqnarray*}
Secondly,
\begin{eqnarray*}
-\ce((e^{-\psi}-1)\varphi,\varphi)
    -\ce(\varphi,(e^{\psi}-1)\varphi)
& = & \int d\sigma^{(\psi,\varphi)}_{1,2}(x_1,x_2)\,(e^{-x_1}-e^{x_1})x_2\\*[5pt]
&&\hspace{5mm} {}-\int d\sigma^{(\psi,\varphi)}_{2,2}(x_1,x_2)\,((e^{-x_1}-1)+(e^{x_1}-1))
\;\;\;.
\end{eqnarray*}
Therefore, by addition,
\[
\ce(\varphi,\varphi)-\ce(e^{-\psi}\varphi,e^{\psi}\varphi)
=\int d\sigma^{(\psi,\varphi)}_{1,1}(x_1,x_2)\,x_2^2
\;\;\;.
\]
Thirdly, two more applications of Proposition~\ref{pdirlocal6} give
\[
\ce(\psi\varphi^2,\psi)
= \int d\sigma^{(\psi,\varphi)}_{1,1}(x_1,x_2)\,x_2^2
   + 2\int d\sigma^{(\psi,\varphi)}_{1,2}(x_1,x_2)\,x_1x_2
\]
and 
\[
2^{-1}\ce(\psi^2,\varphi^2)
=2\int d\sigma^{(\psi,\varphi)}_{1,2}(x_1,x_2) \, x_1 \, x_2
\;\;\;.
\]
Therefore, by subtraction, 
\[
\ci_\psi(\varphi^2)
= \ce(\psi\varphi^2,\psi)-2^{-1}\ce(\psi^2,\varphi^2)
= \int d\sigma^{(\psi,\varphi)}_{1,1}(x_1,x_2)\,x_2^2
\]
which gives the  identity in the proposition.

Finally suppose $\varphi \geq 0$.
Then one calculates similarly that 
\[
\ci_\psi(\varphi)
= \int d\sigma^{(\psi,\varphi)}_{1,1}(x_1,x_2)\,x_2
\;\;\;.  \]
But $\supp \sigma^{(\psi,\varphi)}_{1,1} \subseteq [-\|\psi\|_\infty,\|\psi\|_\infty] \times [0,\|\varphi\|_\infty]$  
by Proposition~\ref{pdirlocal6}.
Hence 
\[
\ci_\psi(\varphi)
= \int d\sigma^{(\psi,\varphi)}_{1,1}(x_1,x_2)\, (x_2 \vee 0)
\;\;\;.  \]
Then the inequality follows immediately from the last part 
of Proposition~\ref{pdirlocal6}.\hfill$\Box$

\ruimte

The following lemma is useful.

\begin{lemma} \label{ldfd208}
Let $\ce$ be a  local Dirichlet form.
\begin{tabel}
\item \label{ldfd208-1}
If $\psi_1,\psi_2,\varphi \in D(\ce) \cap L_\infty(X\!:\!\Ri)$ with $\varphi \geq 0$
then
\[
\ci_{\psi_1 + \psi_2}(\varphi)^{1/2}
\leq \ci_{\psi_1}(\varphi)^{1/2} + \ci_{\psi_2}(\varphi)^{1/2} 
\;\;\; .  \]
\item \label{ldfd208-2}
If $\chi,\psi,\varphi \in D(\ce) \cap L_\infty(X\!:\!\Ri)$ then
\[
|\ci_{\chi \, \psi}(\varphi)|^{1/2}
\leq  \ci_\psi(\chi^2 \, |\varphi|)^{1/2}
   +  \ci_\chi(\psi^2 \, |\varphi|)^{1/2}
\;\;\; .  \]
\end{tabel}
\end{lemma}
\proof\
For all $\psi_1,\psi_2,\varphi \in D(\ce) \cap L_\infty(X\!:\!\Ri)$ set 
\[
\ci_{\psi_1,\psi_2}(\varphi)
= 2^{-1} \ce(\psi_1 \, \varphi,\psi_2) + 2^{-1} \ce(\psi_2 \, \varphi,\psi_1)
   - 2^{-1} \ce(\psi_1 \, \psi_2, \varphi)
\;\;\; .  \]
Then
\[
\Big( (\psi_1,\varphi_1) , (\psi_2,\varphi_2) \Big)
\mapsto \ci_{\psi_1,\psi_2}(\varphi_1 \, \varphi_2)
\]
is a positive symmetric bilinear form on $(D(\ce) \cap L_\infty(X\!:\!\Ri))^2$
by \cite{BH} Proposition~I.4.1.1.
Hence Statement~\ref{ldfd208-1} follows by the corresponding norm triangle inequality
applied to the vectors $(\psi_1,\varphi^{1/2})$ and $(\psi_2,\varphi^{1/2})$.

It follows as in the proof of Proposition~\ref{pdirlocal5} that 
the bilinear form satisfies the Leibniz rule
\[
\ci_{\psi_1 \psi_3, \psi_2}(\varphi)
= \ci_{\psi_1 , \psi_2}(\psi_3 \, \varphi)
   + \ci_{\psi_3 , \psi_2}(\psi_1 \, \varphi)
\]
for all $\psi_1,\psi_2,\psi_3,\varphi \in D(\ce) \cap L_\infty(X\!:\!\Ri)$.
But the Cauchy--Schwarz inequality states that 
\[
|\ci_{\psi_1,\psi_2}(\varphi_1 \, \varphi_2)|
\leq \ci_{\psi_1}(\varphi_1^2)^{1/2} \, \ci_{\psi_2}(\varphi_2^2)^{1/2}
\]
for all $\psi_1,\psi_2,\varphi_1,\varphi_2 \in D(\ce) \cap L_\infty(X\!:\!\Ri)$.
Now let $\chi,\psi,\varphi \in D(\ce) \cap L_\infty(X\!:\!\Ri)$.
Then 
\begin{eqnarray*}
|\ci_{\chi \, \psi}(\varphi)|
\leq \ci_{\chi \, \psi}(|\varphi|)  
& = & \ci_\chi(\psi^2 \, |\varphi|) + 2 \,\ci_{\chi,\psi}(\chi \, \psi \, |\varphi|)
   + \ci_\psi(\chi^2 \, |\varphi|)  \\[5pt]
& \leq & \Big(  \ci_\psi(\chi^2 \, |\varphi|)^{1/2}
   +  \ci_\chi(\psi^2 \, |\varphi|)^{1/2} \Big)^2
\end{eqnarray*}
and Statement~\ref{ldfd208-2} follows.
\hfill$\Box$

\ruimte

We assume from now on in this and the next section that the  Dirichlet form $\ce$ is 
local,  $\one \in D(\ce)_{\rm loc}$ and 
$D(\ce) \cap L_{\infty,c}(X)$ is a core for $\ce$.
These assumptions are valid if 
$X$ is a locally compact separable
metric space and $\mu$ is a Radon measure such that  $\supp \mu = X$,
the Dirichlet form is regular and 
\cite{FOT}-strongly local.

\begin{lemma} \label{ldfd3n08}
If $\varphi \in D(\ce)$ and $\psi \in D(\ce)_{\rm loc} \cap L_\infty(X\!:\!\Ri)$
with $|||\widehat \ci_\psi||| < \infty$ then 
$\psi \, \varphi \in D(\ce)$ and 
\[
\ce(\psi \, \varphi)^{1/2}
\leq |||\widehat \ci_\psi|||^{1/2} \, \|\varphi\|_2
   + \|\psi\|_\infty \, \ce(\varphi)^{1/2}
\;\;\; .  \]
\end{lemma}
\proof\
First assume that $\varphi \in D(\ce) \cap L_{\infty,c}(X\!:\!\Ri)$.
Since $\one \in D(\ce)_{\rm loc}$ there exists a $\chi \in D(\ce) \cap L_\infty(X\!:\!\Ri)$ 
such that $0 \leq \chi \leq 1$ and 
$\chi|_{\supp \varphi} = 1$.
Moreover, there is a $\hat \psi \in D(\ce)$ such that 
$\psi|_{\supp \varphi} = \hat \psi|_{\supp \varphi}$.
We may assume that $\hat \psi \in L_\infty(X\!:\!\Ri)$ and $\|\hat \psi\|_\infty \leq \|\psi\|_\infty$.
Then it follows from  locality and Lemma~\ref{ldfd208}.\ref{ldfd208-2} that 
\begin{eqnarray*}
\ce(\psi \, \varphi)^{1/2}
& = & \ce(\hat \psi \, \varphi)^{1/2}
= \ci_{\hat \psi \, \varphi}(\chi)^{1/2}  \\[5pt]
& \leq & \ci_{\hat \psi}(\varphi^2 \, \chi)^{1/2} + \ci_\varphi(\hat \psi^2 \, \chi)^{1/2}
= \widehat \ci_\psi(\varphi^2 \, \chi)^{1/2} + \ci_\varphi(\hat \psi^2 \, \chi)^{1/2}  \\[5pt]
& \leq & \widehat \ci_\psi(\varphi^2)^{1/2} + \|\hat \psi^2 \, \chi\|_\infty^{1/2} \, \ce(\varphi)^{1/2}
\leq |||\widehat \ci_\psi|||^{1/2} \, \|\varphi\|_2 + \|\psi\|_\infty \, \ce(\varphi)^{1/2}
\;\;\; .  
\end{eqnarray*}
Now let $\varphi \in D(\ce)$.
There exists a sequence $\varphi_1,\varphi_2,\ldots \in D(\ce) \cap L_{\infty,c}(X\!:\!\Ri)$ such that 
$\lim_{n \to \infty} \|\varphi - \varphi_n\|_2 = 0$ and 
$\lim_{n \to \infty} \ce(\varphi - \varphi_n) = 0$.
Then $\lim_{n \to \infty} \psi \, \varphi_n = \psi \, \varphi$ in $L_2(X)$.
Moreover, it follows from the above estimates that 
$n \mapsto \psi \, \varphi_n$ is a Cauchy sequence in $D(\ce)$.
Since $\ce$ is closed one deduces that 
$\psi \, \varphi \in D(\ce)$ and the lemma is established.\hfill$\Box$

\begin{cor} \label{cdfd3n09}
If $\varphi \in D(\ce)$ and $\psi \in D(\ce)_{\rm loc} \cap L_\infty(X\!:\!\Ri)$
with $|||\widehat \ci_\psi||| < \infty$ then 
$e^\psi \varphi \in D(\ce)$ and 
\[
\ce(e^\psi \varphi)^{1/2}
\leq e^{\|\psi\|_\infty} \Big( |||\widehat \ci_\psi|||^{1/2} \, \|\varphi\|_2
                                 +  \ce(\varphi)^{1/2} \Big)
\;\;\; .  \]
\end{cor}
\proof\
It follows by induction from Lemma~\ref{ldfd3n08} that 
\[
\ce(\psi^n \, \varphi)^{1/2}
\leq n \, |||\widehat \ci_\psi|||^{1/2} \, \|\psi\|_\infty^{n-1} \, \|\varphi\|_2
   + \|\psi\|_\infty^n \, \ce(\varphi)^{1/2}
\]
for all $n \in \Ni$, $\varphi \in D(\ce)$ and $\psi \in D(\ce)_{\rm loc} \cap L_\infty(X\!:\!\Ri)$
with $|||\widehat \ci_\psi||| < \infty$.
Since $\ce(\tau)^{1/2} = \|H^{1/2} \tau\|_2$ for all $\tau \in D(\ce) = D(H^{1/2})$ and $H^{1/2}$ is self-adjoint
it follows that $e^\psi \varphi \in D(\ce)$ and 
\begin{eqnarray*}
\ce(e^\psi \varphi)^{1/2}
& \leq & \sum_{n=0}^\infty n!^{-1} \Big( n \, |||\widehat \ci_\psi|||^{1/2} \, \|\psi\|_\infty^{n-1} \, \|\varphi\|_2
   + \|\psi\|_\infty^n \, \ce(\varphi)^{1/2} \Big)  \\[5pt]
& = & e^{\|\psi\|_\infty} \Big( |||\widehat \ci_\psi|||^{1/2} \, \|\varphi\|_2
                                 +  \ce(\varphi)^{1/2} \Big)
\end{eqnarray*}
as required.\hfill$\Box$

\begin{lemma} \label{ldfd3n10}
If $\psi \in D(\ce)_{\rm loc} \cap L_\infty(X\!:\!\Ri)$
with $|||\widehat \ci_\psi||| < \infty$ then 
\begin{equation}
- \ce(e^{-\psi} \varphi, e^\psi \varphi)
\leq |||\widehat \ci_\psi||| \, \|\varphi\|_2^2
\label{eldfd3n10;1}
\end{equation}
for all $\varphi \in D(\ce)$.
\end{lemma}
\proof\
It follows by definition together with Proposition~\ref{pdirlocal5} that
\[
\ci_\psi(\varphi^2)
= \ce(\varphi^2\psi,\psi) - 2^{-1} \ce(\varphi^2,\psi^2)
= \ce(\varphi) - \ce(e^{-\psi} \varphi, e^\psi \varphi)
\]
for all $\varphi,\psi\in D(\ce) \cap L_\infty(X\!:\!\Ri)$.
Therefore 
\[
- \ce(e^{-\psi} \varphi, e^\psi \varphi)
= - \ce(\varphi) + \widehat \ci_\psi(\varphi^2)
\leq \widehat \ci_\psi(\varphi^2)
\]
and (\ref{eldfd3n10;1}) is valid
for all $\psi \in D(\ce)_{\rm loc} \cap L_\infty(X\!:\!\Ri)$ and $\varphi \in D(\ce) \cap L_{\infty,c}(X\!:\!\Ri)$.
But if $\psi \in D(\ce)_{\rm loc} \cap L_\infty(X\!:\!\Ri)$ with $|||\widehat \ci_\psi||| < \infty$ then 
the maps $\varphi \mapsto e^\psi \varphi$ and 
$\varphi \mapsto e^{-\psi} \varphi$ are continuous from $D(\ce)$ into $D(\ce)$ by Corollary~\ref{cdfd3n09}.
Moreover, $D(\ce) \cap L_{\infty,c}(X)$ is dense in $D(\ce)$.
Hence the bounds (\ref{eldfd3n10;1}) extend to all $\varphi \in D(\ce)$ and 
$\psi \in D(\ce)_{\rm loc} \cap L_\infty(X\!:\!\Ri)$ with $|||\widehat \ci_\psi||| < \infty$.\hfill$\Box$

\ruimte

Next for all $\psi \in L_\infty(X)$ define the multiplication operator $M_\psi \colon L_2(X) \to L_2(X)$ by
$M_\psi \varphi = e^\psi \varphi$.

\begin{prop} \label{pdavgaf204}
If $\psi \in D(\ce) \cap L_\infty(X\!:\!\Ri)$ and $|||\widehat \ci_\psi||| < \infty$ then 
\[
\|M_\psi \, S_t \, M_\psi^{-1}\|_{2 \to 2}
\leq e^{|||\widehat \ci_\psi||| t}
\]
for all $t > 0$.
\end{prop}
\proof\
Let $t > 0$.
It follows from Lemma~\ref{ldfd3n10} that 
\[
\|e^\psi S_t \varphi\|_2^2 - \|e^\psi  \varphi\|_2^2
= - 2 \int_0^t ds \, \ce(S_s \varphi, e^{2 \psi} S_s \varphi)  
\leq 2 \int_0^t ds \, |||\widehat \ci_\psi||| \, \|e^\psi S_s \varphi\|_2^2
\]
for all $\varphi \in L_2(X)$.
Then it follows from Gronwall's lemma that 
\[
\|e^\psi S_t \varphi\|_2 
\leq e^{|||\widehat \ci_\psi||| t} \|e^\psi \varphi\|_2
\]
for all $t > 0$.\hfill$\Box$

\ruimte 

Since $\psi\mapsto \ci_\psi(\varphi)$ is a quadratic form one has 
$|||\widehat \ci_{\rho \psi}||| = \rho^2 |||\widehat \ci_\psi|||$ for all $\rho \in \Ri$.
It is now easy to complete the proof of the second theorem.

\ruimte

\noindent
{\bf Proof of Theorem~\ref{tdfd3.1}} \hspace{5pt}\
Let $\psi \in D_0(\ce)$.
Then $|||\widehat \ci_\psi||| \leq  1$ and 
\begin{eqnarray*}
|(\varphi_A, S_t \varphi_B)|
& \leq & \|e^{-\rho \psi} \varphi_A\|_2 \, \|M_{\rho \psi} \, S_t \, M_{\rho \psi}^{-1}\|_{2 \to 2} \,
   \|e^{\rho \psi} \varphi_B\|_2  \\[5pt]
& \leq & e^{|||\widehat \ci_{\rho \psi}||| t} \, e^{-d_{\rho \psi}(A;B)} \, 
                 \|\varphi_A\|_2 \, \|\varphi_B\|_2  
 \leq  e^{\rho^2 t} \, e^{- \rho d_\psi(A;B)} \, \|\varphi_A\|_2 \, \|\varphi_B\|_2 
\end{eqnarray*}
for all $\rho > 0$.
Minimizing over $\psi$ gives
\[
|(\varphi_A, S_t \varphi_B)|
\leq  e^{\rho^2 t} \, e^{- \rho d(A;B)} \, \|\varphi_A\|_2 \, \|\varphi_B\|_2 
\]
for all $\rho > 0$.
If $d(A\,;B) = \infty$ then $(\varphi_A, S_t \varphi_B) = 0$ and the theorem follows.
Finally, if $d(A\,;B) < \infty$ choose $\rho = (2t)^{-1} d(A\,;B)$.\hfill$\Box$

\ruimte

One immediate corollary of the theorem is that the corresponding wave equation has
a finite speed of propagation by  the reasoning of \cite{Sik3}.

\begin{cor}\label{cdfd320}
Let $H$ be the positive self-adjoint operator associated with a 
local Dirichlet form $\ce$ on $L_2(X)$ such that 
$\one \in D(\ce)_{\rm loc}$ and 
$D(\ce) \cap L_{\infty,c}(X)$ is a core for~$\ce$.
If $A,B$ are measurable subsets of $X$ then
\[
(\varphi_A,\cos (tH^{1/2})\varphi_B)=0
\]
for all $\varphi_A\in L_2(A)$, $\varphi_B\in L_2(B)$ and $t\in [-d(A\,;B), d(A\,;B)]$.
\end{cor}
\proof\
This follows immediately from Theorem~\ref{tdfd3.1} and Lemma~3.3 of \cite{ERSZ1}.\hfill$\Box$

\ruimte

We have already remarked in the introduction that Theorem~\ref{tdfd3.1} applies directly to the second-order 
viscosity operators.
Alternatively one may  deduce off-diagonal bounds for operators on open subsets of $\Ri^d$  satisfying  Dirichlet 
or Neumann boundary conditions.
As an illustration consider the Laplacian.

\begin{exam}\label{xdfd328} 
Let $X$ be an  open subset of $\Ri^d$. 
Define $\ce$ on $L_2(X)$ by $D(\ce)=W_0^{1,2}(X)$ and
$\ce(\varphi)=\|\nabla\varphi\|_2^2$.
In particular $D(\ce)$  is the closure of $C_c^\infty(X)$ with respect to the norm
$\varphi\mapsto (\|\nabla\varphi\|_2^2+\|\varphi\|_2^2)^{1/2}$.
It follows that $\ce$ is a regular local Dirichlet form and  the 
assumptions of Theorem~\ref{tdfd3.1} are satisfied.
Therefore Theorem~\ref{tdfd3.1} gives off-diagonal bounds.
The self-adjoint operator corresponding to $\ce$ is the Dirichlet Laplacian on $L_2(X)$.
\end{exam}

\begin{exam}\label{xdfd350} 
Let $\Omega$ be an open subset of $\Ri^d$ with $|\partial \Omega| = 0$.
Set $X = \overline{\Omega}$.
Define the form $\ce$ on $L_2(X) = L_2(\Omega)$  by
$D(\ce) = W^{1,2}(\Omega)$ and $\ce(\varphi) = \|\nabla \varphi\|_2^2$.
Then $\ce$ is a local Dirichlet form.
Fix $\tau \in C_c^\infty(\Ri^d)$ with $0 \leq \tau \leq 1$ and 
$\tau|_{B_e(0;1)} = 1$, where $B_e(0\,;1)$ is the Euclidean unit ball.
For all $n \in \Ni$ define $\tau_n \in C_c^\infty(\Ri^d)$ by 
$\tau_n(x) = \tau(n^{-1} x)$.
Then $\tau_n|_\Omega \in D(\ce)$ for all $n \in \Ni$ and therefore $\one \in D(\ce)_{\rm loc}$.
The space $D(\ce) \cap L_\infty(X)$ is dense in $D(\ce)$ by \cite{FOT} Theorem~1.4.2.(iii).
If $\varphi \in D(\ce) \cap L_\infty(X)$ then $\lim_{n \to \infty} \tau_n|_\Omega \, \varphi = \varphi$
in $W^{1,2}(\Omega) = D(\ce)$.
But for all $n \in \Ni$ the support of the function $\tau_n|_\Omega \, \varphi$, viewed 
as an almost everywhere defined function on $X$, is closed in $X$, and therefore also closed in $\Ri^d$.
Hence this support is compact in $\Ri^d$ and then also compact in $X$.
So $\tau_n|_\Omega \, \varphi \in D(\ce) \cap L_{\infty,c}(X)$ and the 
space $D(\ce) \cap L_{\infty,c}(X)$ is dense in $D(\ce)$.
Therefore Theorem~\ref{tdfd3.1} gives off-diagonal bounds.
The self-adjoint operator corresponding to $\ce$ is the Neumann Laplacian on $L_2(\Omega)$.

Note that we do not assume that $\Omega$ has the segment property. 
In general the Dirichlet form $\ce$ is not regular on $X$.
\end{exam}

\section{Separation} \label{Sdfd4}

In this section we give the proof of Theorem~\ref{tdfd401}.
In \cite{ERSZ1} we established that  for 
degenerate elliptic operators phenomena of separation can occur. 
In particular the corresponding semigroup $S$ is reducible, i.e., it has non-trivial invariant
subspaces.
The theorem  shows that such subspaces can be characterized by the 
 set-theoretic  distance $d(\cdot\,;\cdot)$.

For the proof of the implication \ref{tdfd401-0}$\Rightarrow$\ref{tdfd401-1} in Theorem~\ref{tdfd401}
we need a variation of an argument used to prove Theorem~XIII.44
in \cite{RS4} (see also \cite{bSim6}).
The essence of the argument is contained in the following lemma.

\begin{lemma} \label{ldfd4.1}
Let $A,B$ be measurable subsets of $X$ with finite measure.
If $(S_t\one_A,\one_B) = 0$ for one $t>0$ then $(S_t\one_A,\one_B) = 0$ for all $t>0$.
\end{lemma}
\proof\
Set $P = \{ t \in \langle0,\infty\rangle : (S_t \one_A,\one_B) > 0 \} $.
Let $t \in P$ and $\tau > 0$. 
Since $(S_t\one_A,\one_B) > 0$ it follows that the product 
$(S_t\one_A) \cdot (\one_B)$ is not identically zero.
Let $\varphi=S_t\one_A\wedge \one_B$.
Then $\varphi \neq 0$.
Moreover, 
\[
(S_{t+\tau} \one_A, \one_B)
= (S_\tau(S_t\one_A), \one_B)
\geq (S_\tau\varphi, \varphi)
= \|S_{\tau/2}\varphi\|_2^2
>0
\]
because $S$ is Markovian.
Thus $t+\tau \in P$ and $[t,\infty\rangle \subset P$ for all $t \in P$.

Set $N = \langle0,\infty\rangle \backslash P = \{ t \in \langle0,\infty\rangle : (S_t \one_A,\one_B) = 0 \} $.
If $N \neq \emptyset$ then there is an $a > 0$ such that $\langle0,a\rangle \subset N$.
But the map $z \mapsto (S_z\one_A,\one_B)$ is analytic in the open right half-plane.
Hence if $N \neq \emptyset$ then $(S_t\one_A,\one_B) = 0$ for all $t \in \langle0,\infty\rangle$.\hfill$\Box$

\ruimte

Before we prove Theorem~\ref{tdfd401} we need one more lemma.
Recall that we assume throughout this section that  the Dirichlet form $\ce$ is 
local,  $\one \in D(\ce)_{\rm loc}$ and 
$D(\ce) \cap L_{\infty,c}(X)$ is a core for $\ce$.
A normal contraction on $\Ri$ is a function $F \colon \Ri \to \Ri$
such that 
$|F(x) - F(y)| \leq |x-y|$ for all $x,y \in \Ri$ and $F(0) = 0$.

\begin{lemma} \label{ldfd438}
\mbox{}
\begin{tabel}
\item \label{ldfd438-1}
If $\psi \in D(\ce)_{\rm loc} \cap L_\infty(X\!:\!\Ri)$ and $\lambda \in \Ri$ then 
$\psi + \lambda \, \one \in D(\ce)_{\rm loc} \cap L_\infty(X\!:\!\Ri)$ and 
$\widehat \ci_{\psi + \lambda \one} = \widehat \ci_\psi$.
In particular $|||\widehat \ci_{\psi + \lambda \one}||| = |||\widehat \ci_\psi|||$.
\item \label{ldfd438-2}
If $\psi \in D(\ce)_{\rm loc} \cap L_\infty(X\!:\!\Ri)$ and 
$F$ is a normal contraction then $F \circ \psi \in D(\ce)_{\rm loc} \cap L_\infty(X\!:\!\Ri)$ and 
$|||\widehat \ci_{F \circ \psi}||| \leq |||\widehat \ci_\psi|||$.
\item \label{ldfd438-3}
If $A,B$ are measurable subsets of $X$, $M \in [0,d(A\,;B)] \cap \Ri$ and $\varepsilon > 0$
then there is a $\psi \in D_0(\ce)$ such that 
$0 \leq \psi \leq M$, $\psi(b) = 0$ for a.e.\ $b \in B$ and 
$\psi(a) \geq M - \varepsilon$ for a.e.\ $a \in A$.
\end{tabel}
\end{lemma}
\proof\
Let $\varphi \in D(\ce) \cap L_{\infty,c}(X\!:\!\Ri)$.
There exist $\hat \psi,\chi \in D(\ce) \cap L_\infty(X\!:\!\Ri)$ such that 
$\hat \psi|_{\supp \varphi} = \psi|_{\supp \varphi}$
and 
$\chi|_{\supp \varphi} = 1$.
Then the locality of $\ce$ implies that 
\begin{eqnarray*}
\widehat \ci_{\psi + \lambda \one}(\varphi)
& = & \ci_{\hat \psi + \lambda \chi}(\varphi)
= \ce((\hat \psi + \lambda \chi) \varphi, \hat \psi + \lambda \chi) - 2^{-1} \ce(\varphi,(\hat \psi + \lambda \chi)^2) \\[5pt]
& = & \ce(\hat \psi \, \varphi, \hat \psi) + \lambda \, \ce(\varphi,\hat \psi)
   - 2^{-1} \Big( \ce(\varphi,\hat \psi^2) - 2 \ce(\varphi,\lambda \, \hat \psi \, \chi) \Big)  \\[5pt]
& = & \ce(\hat \psi \, \varphi, \hat \psi) - 2^{-1}  \ce(\varphi,\hat \psi^2)
= \ci_{\hat \psi}(\varphi)
= \widehat \ci_\psi(\varphi)
\;\;\; .  
\end{eqnarray*}
This proves Statement~\ref{ldfd438-1}.

If $F$ is a normal contraction on $\Ri$ then 
$\ci_{F \circ \psi}(\varphi) \leq \ci_\psi(\varphi)$ for all 
$\varphi,\psi \in D(\ce) \cap L_\infty(X\!:\!\Ri)$ with $\varphi \geq 0$
by \cite{BH} Proposition~I.4.1.1.
Then Statement~\ref{ldfd438-2} is an easy consequence.

Finally, there exists a $\psi \in D_0(\ce)$ such that $d_\psi(A\,;B) \geq M - \varepsilon$.
We may assume that $|B| > 0$.
Then $0 \vee (\psi + \|\one_B \, \psi\|_\infty) \wedge M$
satisfies the required conditions.\hfill$\Box$

\ruimte

\noindent
{\bf Proof of Theorem~\ref{tdfd401}}\hspace{5pt}\
The implication \ref{tdfd401-0}$\Rightarrow$\ref{tdfd401-1} is now immediate.
Let  $B\subset  A^{\rm c}$, $A' \subset A$ and suppose that $A'$ and $B$ have finite measure.
Then $(S_t\one_{A'},\one_B)=0$ for one $t>0$ by assumption and for all $t>0$ by
Lemma~\ref{ldfd4.1}.
Hence $(S_t \varphi, \psi) = 0$ for all $t > 0$, $\varphi \in L_2(A)$ and $\psi \in L_2(A^{\rm c})$.
Thus $S_tL_2(A)\subseteq L_2(A)$ for all $t>0$.

The converse implication \ref{tdfd401-1}$\Rightarrow$\ref{tdfd401-0} is obvious.

Next we prove the implication \ref{tdfd401-1}$\Rightarrow$\ref{tdfd401-2}.
For all $\varphi \colon X \to \Ci$ set $\widetilde \varphi = \one_A \, \varphi$.
Then $\widetilde \varphi \in D(\ce)$ and 
$\ce(\psi,\widetilde \varphi) = \ce(\widetilde \psi,\varphi) = \ce(\widetilde \psi,\widetilde \varphi)$
for all $\varphi,\psi \in D(\ce)$ by \cite{ERSZ1}, Lemma~6.3.
Hence
$\ci_{\widetilde \psi}(\varphi) = \ce(\widetilde \psi \, \varphi, \widetilde \psi) - 2^{-1} \ce(\varphi, \widetilde \psi^2)
   = \ci_\psi(\widetilde \varphi)$
for all $\varphi,\psi \in D(\ce) \cap L_\infty(X\!:\!\Ri)$.
Therefore, if $\psi \in D(\ce)_{\rm loc} \cap L_\infty(X\!:\!\Ri)$ and 
$\varphi \in D(\ce) \cap L_{\infty,c}(X\!:\!\Ri)$ then 
$\tilde \psi \in D(\ce)_{\rm loc} \cap L_\infty(X\!:\!\Ri)$
and $\widehat \ci_{\tilde \psi}(\varphi) = \widehat \ci_\psi(\tilde \varphi)$.
But $\widehat \ci_{\one} = 0$ by Lemma~\ref{ldfd438}.\ref{ldfd438-1} since $\ce$ is  local.
So $\widehat \ci_{\one_A}(\varphi) = \widehat \ci_{\tilde{\one}}(\varphi)
= \widehat \ci_{\one}(\tilde \varphi) = 0$ for all 
$\varphi \in D(\ce) \cap L_{\infty,c}(X\!:\!\Ri)$.
So $|||\widehat \ci_{\one_A}||| = 0$.
Then $d(A\,;A^{\rm c}) = \infty$.

The implication \ref{tdfd401-2}$\Rightarrow$\ref{tdfd401-4} is trivial.

Finally suppose that \ref{tdfd401-4} is valid.
Let $\delta \in \langle0,\infty\rangle$ such that $2 \delta < d(A\,;A^{\rm c})$.
By Lemma~\ref{ldfd438}.\ref{ldfd438-3} there exists a 
$\psi \in D(\ce)_{\rm loc} \cap L_\infty(X\!:\!\Ri)$ such that 
$\psi(b) = 0$ for a.e.\ $b \in A^{\rm c}$ and
$\psi(a) \geq \delta$ for a.e.\ $a \in A$.
Then $\one_A = 0 \vee (\delta^{-1} \psi_2) \wedge 1$ and 
since $x \mapsto 0 \vee x \wedge 1$ is a normal contraction it 
follows from Lemma~\ref{ldfd438}.\ref{ldfd438-2} that $\one_A \in D(\ce)_{\rm loc} \cap L_\infty(X\!:\!\Ri)$
and $|||\widehat \ci_{\one_A}||| \leq \delta^{-2} < \infty$.

Now let $\varphi \in D(\ce)$.
It follows from Lemma~\ref{ldfd3n08} that $\one_A \, \varphi \in D(\ce)$.
Then also $\one_{A^{\rm c}} \, \varphi \in D(\ce)$.
For all $t > 0$ one has by Theorem~\ref{tdfd3.1} that 
\begin{eqnarray*}
t^{-1} \, |(\one_A \, \varphi, (I - S_t) (\one_{A^{\rm c}} \, \varphi))|
& = & t^{-1} \, |(\one_A \, \varphi, S_t (\one_{A^{\rm c}} \, \varphi))|  \\[5pt]
& \leq & t^{-1} \, e^{-\delta^2 t^{-1}} \|\one_A \, \varphi\|_2 \, \|\one_{A^{\rm c}} \, \varphi\|_2
\;\;\; .  
\end{eqnarray*}
Hence
\[
\ce(\one_A \, \varphi, \one_{A^{\rm c}} \, \varphi)
= \lim_{t \downarrow 0} 
   t^{-1} \, (\one_A \, \varphi, (I - S_t) (\one_{A^{\rm c}} \, \varphi))
= 0
\]
and 
\[
\ce(\varphi)
= \ce(\one_A \, \varphi) + \ce(\one_{A^{\rm c}} \, \varphi)
\;\;\; .  \]
Then \cite{ERSZ1}, Lemma~6.3, implies that \ref{tdfd401-1} is valid.\hfill$\Box$

\ruimte

Note that the above proof shows that the equivalent statements of 
Theorem~\ref{tdfd401} are also equivalent with the 
statement $\one_A \in D(\ce)_{\rm loc} \cap L_\infty(X\!:\!\Ri)$
and $|||\widehat \ci_{\one_A}||| < \infty$.

\begin{exam} \label{xdfd413}
Let $\delta \in [0,1\rangle$ and
consider the viscosity form $h$ on $\Ri^d$ with $d = 1$ and coefficient
$c_{11} = c_\delta$ where 
\[
c_\delta(x)
= \Big( \frac{x^2}{1+x^2} \Big)^\delta
\;\;\; .  \]
Then $h(\varphi) = \int |\varphi'(x)|^2 \, c_\delta(x)$ for all 
$\varphi\in C_c^\infty(\Ri)$, 
the form $h$ is closable and its closure is a Dirichlet form.
If $S$ is the semigroup generated by the operator associated with the closure and if
 $\delta\in[1/2,1\rangle$ then $S_tL_2(-\infty,0)\subseteq L_2(-\infty,0)$ and 
$S_tL_2(0,\infty)\subseteq L_2(0,\infty)$ for all $t>0$
by \cite{ERSZ1}, Corollary~2.4 and Proposition~6.5.
The assumptions  of Theorems~\ref{tdfd3.1} and \ref{tdfd401} are satisfied and $d(A\,;B)=\infty$
for each pair of measurable subsets $A\subseteq \langle -\infty,0\rangle$ and 
$B \subseteq \langle0,\infty\rangle$.
Note, however, that the corresponding Riemannian distance
\[
d(x\,;y)=\Big|\int^x_yds\,c_\delta(s)^{-1/2}\,\Big|
\]
is finite for all $x,y\in\Ri$.
Therefore the Riemannian distance does not reflect the behaviour of the semigroup.
\end{exam}

\section{Distances}\label{Sdfd5}

In this section we 
derive various general properties of the  distance $d^{(\ce)}(\cdot\,;\cdot)$
 and give several  examples.

Let $\ce$ be a local Dirichlet form on $L_2(X)$.
If $A_1 \subseteq A_2$ and $B$ are measurable then $d_\psi(A_1\,;B) \geq d_\psi(A_2\,;B)$
for all $\psi \in L_\infty(X)$.
Hence $d(A_1\,;B) \geq d(A_2\,;B)$.
Also $d_\psi(B\,;A) = d_{-\psi}(A\,;B)$ for all $\psi \in L_\infty(X)$.
So $d(A\,;B) = d(B\,;A)$ and 
$d(A\,;B_1) \geq d(A\,;B_2)$ whenever $B_1 \subseteq B_2$ are measurable.
Next we consider monotonicity of the distance as a function of the form.

If $h_1$ and $h_2$ are two strongly elliptic forms on $\Ri^d$ and $h_1\geq h_2$ then
the corresponding matrices of coefficients satisfy $C^{(1)}\geq C^{(2)}$.
Therefore $(C^{(1)})^{-1}\leq (C^{(2)})^{-1}$.
Hence the corresponding Riemannian distances satisfy $d_1(x\,;y)\leq d_2(x\,;y)$ 
for all $x,y\in\Ri^d$.
Thus the order of the forms gives the inverse order for the distances.
One can establish a similar result for general Dirichlet forms and the set-theoretic 
distance under subsidiary regularity conditions.

\begin{prop}\label{pdirlocal5.1}
Let $\ce$ and $\cf$ be local Dirichlet forms with $\ce \leq \cf$.
Assume $\one \in D(\cf)_{\rm loc}$, the space 
$D(\cf) \cap L_{\infty,c}(X)$ is a core for  $D(\ce)$ and 
for each  compact $K \subset X$ there is a $\chi \in D(\cf) \cap L_{\infty,c}(X\!:\!\Ri)$
such that $0 \leq \chi \leq 1$, $\chi|_K = 1$ and $|||\widehat I^{(\cf)}_\chi||| < \infty$.
Then 
\[
d^{(\cf)}(A\,;B)\leq d^{(\ce)}(A\,;B)
\]
for all measurable  $A,B\subseteq X$.
\end{prop}
\proof\
It suffices to prove that   $|||\widehat \ci^{(\ce)}_\psi|||\leq|||\widehat \ci^{(\cf)}_\psi|||$
for each $\psi \in D(\cf)_{\rm loc} \cap L_\infty(X\!:\!\Ri)$ 
because the statement of the proposition then follows from the definition of the distance.

If $\psi,\varphi \in D(\cf) \cap L_\infty(X\!:\!\Ri)$ with $\varphi \geq 0$ then 
$\ci^{(\ce)}_\psi(\varphi)\leq \ci^{(\cf)}_\psi(\varphi)$
by Proposition~\ref{pdirlocal5}.
Therefore $\widehat \ci^{(\ce)}_\psi(\varphi) \leq \widehat \ci^{(\cf)}_\psi(\varphi)$
for all $\psi \in D(\cf)_{\rm loc} \cap L_\infty(X\!:\!\Ri)$ and 
$\varphi \in D(\cf) \cap L_{\infty,c}(X\!:\!\Ri)$ with $\varphi \geq 0$.
Fix $\psi \in D(\cf)_{\rm loc} \cap L_\infty(X\!:\!\Ri)$.

If $|||\widehat \ci^{(\cf)}_\psi|||=\infty$ there is nothing to prove, so we may assume that 
$|||\widehat \ci^{(\cf)}_\psi||| < \infty$.

It follows as in the proof of Lemma~\ref{ldfd3n08} that 
\begin{eqnarray*}
\ce(\psi \varphi)^{1/2}
& \leq & \widehat \ci^{(\ce)}_\psi(\varphi^2)^{1/2} + \|\psi\|_\infty \, \ce(\varphi)^{1/2}  \\[5pt]
& \leq & \widehat \ci^{(\cf)}_\psi(\varphi^2)^{1/2} + \|\psi\|_\infty \, \ce(\varphi)^{1/2}  
 \leq  |||\widehat \ci^{(\cf)}_\psi||| \, \|\varphi\|_2 + \|\psi\|_\infty \, \ce(\varphi)^{1/2} 
\end{eqnarray*}
for all $\varphi \in D(\cf) \cap L_{\infty,c}(X\!:\!\Ri)$.
Since $\ce$ is closed and $D(\cf) \cap L_{\infty,c}(X)$ is a core
it then follows that $\psi \varphi \in D(\ce)$ for all $\varphi \in D(\ce)$
and 
\[
\ce(\psi \varphi)^{1/2}
\leq |||\widehat \ci^{(\cf)}_\psi||| \, \|\varphi\|_2 + \|\psi\|_\infty \, \ce(\varphi)^{1/2} 
\]
for all $\varphi \in D(\ce)$.

Let $\varphi_0 \in D(\ce) \cap L_{\infty,c}(X\!:\!\Ri)$.
By assumption there exists a $\chi \in D(\cf) \cap L_{\infty,c}(X\!:\!\Ri)$
such that $0 \leq \chi \leq 1$, $\chi|_{\supp \varphi_0} = 1$ and $|||\widehat \ci^{(\cf)}_\chi||| < \infty$.
Then by the above argument with $\psi$ replaced by $\chi$ one deduces that
\[
\ce(\chi \varphi)^{1/2}
\leq |||\widehat \ci^{(\cf)}_\chi||| \, \|\varphi\|_2 + \ce(\varphi)^{1/2} 
\]
for all $\varphi \in D(\ce)$.
There exists a $\hat \psi \in D(\ce) \cap L_\infty(X\!:\!\Ri)$ such that 
$\hat \psi|_{\supp \chi} = \psi|_{\supp \chi}$.
Then there is a $c > 0$ such that 
\begin{eqnarray*}
|\widehat \ci^{(\ce)}_\psi(\chi \, \varphi)|
& = & |\ci^{(\ce)}_{\hat \psi}(\chi \, \varphi)|
\leq \ce(\hat \psi)^{1/2} \ce(\hat \psi \, \chi \, \varphi)^{1/2}
   +2^{-1} \ce(\hat \psi^2)^{1/2} \, \ce(\chi \, \varphi)^{1/2}   \nonumber  \\[5pt]
& = & \ce(\hat \psi)^{1/2} \ce(\psi \, \chi \, \varphi)^{1/2}
   +2^{-1} \ce(\hat \psi^2)^{1/2} \, \ce(\chi \, \varphi)^{1/2}   \nonumber  \\[5pt]
&\leq & \ce(\hat \psi)^{1/2} \Big( |||\widehat \ci^{(\cf)}_\psi||| \, \|\chi \, \varphi\|_2 + \|\psi\|_\infty \, \ce(\chi \, \varphi)^{1/2}  \Big)
   + 2^{-1}\ce(\psi^2)^{1/2} \, \ce(\chi \, \varphi)^{1/2}  \\[5pt]
& \leq & c \Big( \ce(\varphi)^{1/2} + \|\varphi\|_2 \Big)
\end{eqnarray*}
uniformly for all $\varphi \in D(\ce) \cap L_{\infty,c}(X\!:\!\Ri)$.
Since the space $D(\cf) \cap L_{\infty,c}(X\!:\!\Ri)$ is dense in $D(\ce)$
there are $\varphi_1,\varphi_2,\ldots \in D(\cf) \cap L_{\infty,c}(X\!:\!\Ri)$ such that 
$\lim_{n \to \infty} (\|\varphi_0 - \varphi_n\|_2 + \ce(\varphi_0 - \varphi_n)) = 0$.
Then 
\[
\lim_{n \to \infty} \widehat \ci^{(\ce)}_\psi(\chi \varphi_n) 
= \widehat \ci^{(\ce)}_\psi(\chi \varphi_0) 
= \widehat \ci^{(\ce)}_\psi(\varphi_0) 
\;\;\; .  \]
On the other hand, 
\[
|\widehat \ci^{(\ce)}_\psi(\chi \varphi_n)|
\leq |||\widehat \ci^{(\cf)}_\psi||| \, \|\chi \varphi_n\|_1
\leq |||\widehat \ci^{(\cf)}_\psi||| \, \|\varphi_0\|_1 + |||\widehat \ci^{(\cf)}_\psi||| \, \|\chi\|_2 \, \|\varphi_0 - \varphi_n\|_2
\]
for all $n \in \Ni$.
Taking the limit $n \to \infty$ gives $|\widehat \ci^{(\ce)}_\psi(\varphi_0)| \leq |||\widehat \ci^{(\cf)}_\psi||| \, \|\varphi_0\|_1$.
The required monotonicity of the norms is immediate.\hfill$\Box$

\ruimte
One can  apply the proposition to a general elliptic form $h$ as in (\ref{edfd1.1})
to obtain a lower bound on 
the distance.
Then the viscosity form $h_0$ satisfies 
 $h_0\leq h\leq \|C\|\,l$  where $l$ is the form associated with the Laplacian.
Therefore 
\[
d^{(h_0)}(A\,;B)\geq \|C\|^{-1/2}d^{(l)}(A\,;B)
\]
for all measurable $A,B\subseteq \Ri^d$.
Alternatively if $h$ is a  strongly  elliptic form then there is a $\mu > 0$ such that $C\geq\mu I$
and $\mu\,l\leq h\leq \|C\|\,l$.
It then follows that 
$d^{(h)}(\cdot\,;\cdot)$ is equivalent to the  distance $d^{(l)}(\cdot\,;\cdot)$.
Specifically,
 \[
 \|C\|^{-1/2}d^{(l)}(A\,;B)\leq d^{(h)}(A\,;B)\leq \mu^{-1/2}d^{(l)}(A\,;B)
  \]
for all measurable sets $A$ and $B$.

In the case of the Laplacian one can explicitly identify the distance.

\begin{exam} \label{xdfd318}
If $A$ and $B$ are non-empty open subsets of $\Ri^d$ then
\begin{equation}
d^{(l)}(A\,;B)=\inf_{x\in A}\,\inf_{y\in B}|x-y|
\label{edfd5.99}
\;\;\;.
\end{equation}
First, observe that if 
$\varphi,\psi\in  D(l)\cap L_\infty(\Ri^d\!:\!\Ri) = W^{1,2}\cap L_\infty(\Ri^d\!:\!\Ri)$
then 
\[
\ci^{(l)}_\psi(\varphi)
=(\nabla(\varphi\psi),\nabla \psi)-2^{-1}(\nabla \varphi,\nabla(\psi^2))
=\int_{\Ri^d}dx\,\varphi(x) \, |(\nabla\psi)(x)|^2
\;\;\;.
\]
Secondly, let $\psi\in W^{1,2}_{\rm loc}\cap L_\infty(\Ri^d\!:\!\Ri)$ and 
$K \subset \Ri^d$ compact.
Choose $\hat\psi\in W^{1,2}\cap L_\infty$
such that $\psi|_K=\hat\psi|_K$.
If $|||\widehat \ci^{(l)}_\psi|||\leq 1$
then 
\[
\Big| \int_K dx\,\varphi(x) \, |(\nabla\psi)(x)|^2 \Big|
= \Big| \int_K dx\,\varphi(x) \, |(\nabla{\hat \psi})(x)|^2 \Big|
= | \ci^{(l)}_{\hat \psi}(\varphi)|
= |\widehat \ci^{(l)}_\psi(\varphi)|
\leq \|\varphi\|_1
\]
for all $\varphi\in W^{1,2}\cap L_\infty(\Ri^d\!:\!\Ri) \cap L_1$ with  $\supp \varphi \subseteq K$.
Therefore 
$\sup_{x\in K}|(\nabla\psi)(x)|\leq 1$ uniformly for all $K$
and   $\psi\in W^{1,\infty}$.
Thus $D_0(l)=\{\psi\in W^{1,\infty}(\Ri^d\!:\!\Ri) : \|\nabla\psi\|_\infty\leq 1\}$.
But it is well known that 
\[
|x-y|=\sup\{|\psi(x)-\psi(y)|:W^{1,\infty}(\Ri^d\!:\!\Ri)\, , \,\|\nabla\psi\|_\infty\leq 1\}
\]
(see, for example, \cite{JSC1}, Proposition~3.1).
Then   (\ref{edfd5.99}) follows immediately.
\end{exam}

The next proposition compares  forms under a rather different regularity assumption.
We define $D(\cf)$ to be an {\bf ideal} of $D(\ce)$ if $D(\cf) \subseteq D(\ce)$ and 
\[
\Big(D(\cf)_{\rm loc}\cap L_\infty(X\!:\!\Ri)\Big)\,
\Big(D(\ce)\cap L_{\infty;c}(X\!:\!\Ri)\Big)
\subseteq 
\Big(D(\cf)_{\rm loc}\cap L_\infty(X\!:\!\Ri)\Big)
\;\;\;.
\]
Note that if $\one\in D(\cf)_{\rm loc}$ 
and if $D(\cf)$ is an ideal of $D(\ce)$ then 
\[
\psi=\one\,\psi \in \Big(D(\cf)_{\rm loc}\cap L_\infty(X\!:\!\Ri)\Big)\,
\Big(D(\ce)\cap L_{\infty;c}(X\!:\!\Ri)\Big)
\subseteq 
\Big(D(\cf)_{\rm loc}\cap L_{\infty;c}(X\!:\!\Ri)\Big)
\]
for all $\psi\in D(\ce)\cap L_{\infty;c}(X\!:\!\Ri)$.
So $D(\cf)_{\rm loc}\cap L_{\infty;c}(X\!:\!\Ri)=D(\ce)_{\rm loc}\cap L_{\infty;c}(X\!:\!\Ri)$.

\begin{prop}\label{pdirlocal5.2}
Let  $\ce$ and $\cf$ be local Dirichlet forms such that 
$\ce \leq \cf$.
Assume $\one\in D(\cf)_{\rm loc}$, the space $D(\cf)$ is an ideal of $D(\ce)$
and $D(\cf) \cap L_{\infty,c}(X)$ is dense in $L_1$.
Then
\[
d^{(\cf)}(A\,;B)\leq d^{(\ce)}(A\,;B)
\]
for all measurable  $A,B\subset X$.\end{prop}
\proof\
It follows from the definition of the distance and Lemma~\ref{ldfd438} that it suffices 
 to prove that   $|||\widehat \ci^{(\ce)}_\psi|||\leq|||\widehat \ci^{(\cf)}_\psi|||$
for all positive $\psi \in D(\cf)_{\rm loc} \cap L_\infty(X\!:\!\Ri)$.

First, 
one has again  $\widehat \ci^{(\ce)}_\psi(\varphi)\leq \widehat \ci^{(\cf)}_\psi(\varphi)$
for all $\psi \in D(\cf)_{\rm loc} \cap L_\infty(X\!:\!\Ri)$ and 
$\varphi \in D(\cf) \cap L_{\infty,c}(X\!:\!\Ri)$ with $\varphi\geq0$.

Fix a positive  $\psi \in D(\cf)_{\rm loc} \cap L_\infty(X\!:\!\Ri)$.
If $|||\widehat \ci^{(\cf)}_\psi|||=\infty$ there is nothing to prove.
So we may assume $|||\widehat \ci^{(\cf)}_\psi|||<\infty$.

Since $\one\in D(\cf)_{\rm loc}$ it follows that $\one\in D(\ce)_{\rm loc}$.
It then follows  from Lemma~\ref{ldfd438} that 
\[
\widehat \ci^{(\ce)}_\psi(\varphi)=\widehat \ci^{(\ce)}_{\one+\psi}(\varphi)
= 4\,\widehat \ci^{(\ce)}_{(\one+\psi)^{1/2}}((\one+\psi)\varphi)
\]
for all $\varphi\in D(\ce)\cap L_{\infty;c}(X\!:\!\Ri)$
 by the Leibniz rule.
But $(\one+\psi)\varphi\in D(\cf)_{\rm loc}\cap L_{\infty;c}(X\!:\!\Ri)$
by the ideal property.
Therefore
\[
\widehat \ci^{(\ce)}_\psi(\varphi)
\leq 
 4\,\widehat \ci^{(\cf)}_{(\one+\psi)^{1/2}}((\one+\psi)\varphi)
\]
for all positive $\varphi \in D(\ce) \cap L_{\infty,c}(X\!:\!\Ri)$.
Hence 
\[
|||\widehat \ci^{(\ce)}_\psi|||\leq 4\,(1+\|\psi\|_\infty)\,
|||\widehat \ci^{(\cf)}_{(\one+\psi)^{1/2}}|||
\leq 4\,(1+\|\psi\|_\infty)\,
|||\widehat \ci^{(\cf)}_\psi|||
\]
where the last inequality follows from Lemma~\ref{ldfd438}.\ref{ldfd438-2}.
In particular $|||\widehat \ci^{(\ce)}_\psi|||<\infty$.

Finally, since 
$|||\widehat \ci^{(\ce)}_\psi|||<\infty$ one has 
\[
|\widehat \ci^{(\ce)}_\psi(\varphi)|\leq \widehat \ci^{(\ce)}_\psi(|\varphi|)
\leq |||\widehat \ci^{(\ce)}_\psi|||\,\|\varphi\|_1
\]
for all  $\varphi \in D(\ce) \cap L_{\infty,c}(X\!:\!\Ri)$.
But then there is a $\Gamma \in L_\infty$ such that $\|\Gamma\|_\infty = |||\widehat \ci^{(\ce)}_\psi|||$ and
\[
\widehat \ci^{(\ce)}_\psi(\varphi)
= \int \Gamma \, \varphi
\]
for all $\varphi \in D(\ce) \cap L_{\infty,c}(X\!:\!\Ri)$.
Therefore 
\[
\Big| \int \Gamma \, \varphi \,\Big|
= |\widehat \ci^{(\ce)}_\psi(\varphi)|
\leq |\widehat \ci^{(\ce)}_\psi(|\varphi|)|
\leq |\widehat \ci^{(\cf)}_\psi(|\varphi|)|
\leq |||\widehat \ci^{(\cf)}_\psi||| \, \|\varphi\|_1
\]
for all $\varphi \in D(\cf) \cap L_{\infty,c}(X\!:\!\Ri)$.
Since the latter space is dense in $L_1(X\!:\!\Ri)$ one deduces that 
$\|\Gamma\|_\infty \leq |||\widehat \ci^{(\cf)}_\psi|||$.
Hence $|||\widehat \ci^{(\ce)}_\psi||| \leq |||\widehat \ci^{(\cf)}_\psi|||$.
The required monotonicity of the distances is immediate.\hfill$\Box$

\ruimte

The proposition allows comparison of elliptic operators with different boundary 
conditions.

\begin{exam}\label{xdfd517}
Let $X$ be an open subset of $\Ri^d$.
Define the forms of the Neumann and Dirichlet Laplacians on $L_2(X)$ by
$h_N(\varphi)=\|\nabla\varphi\|_2^2$ with $D(h_N)=W^{1,2}(X)$ and 
$h_D(\varphi)=\|\nabla\varphi\|_2^2$ with $D(h_D)=W^{1,2}_0(X)$.
Then $h_N\leq h_D$, $\one\in D(h_D)_{\rm loc}$, $D(h_D)$ is an ideal of $D(h_N)$
and $D(h_D) \cap L_{\infty,c}(X)$ is dense in $L_1(X)$.
Therefore one deduces from Proposition~\ref{pdirlocal5.2} that $d^{(h_D)}(A\,;B)\leq d^{(h_N)}(A\,;B)$
for all measurable $A,B\subseteq X$.
Note that Proposition~\ref{pdirlocal5.1} does not apply to this example since 
$W^{1,2}_0(X) \cap L_{\infty,c}(X)$ is not a core for $D(h_N)$. 
\end{exam}

The argument used in Example~\ref{xdfd318} allows one to identify the distances associated
with $h_N$ and $h_D$ with the geodesic distance in the Euclidean metric.
In particular the distance is independent of the boundary conditions.

\begin{exam}\label{xdfd518}
Let $X$ be an  open subset of $\Ri^d$ and $A,B$ non-empty open subsets of $X$.
If $X$ is disconnected then 
$d^{(h_N)}(A\,;B)=\infty=d^{(h_D)}(A\,;B)$
whenever $A$ and $B$ are in separate components  by Theorem~\ref{tdfd401}.
Hence we may assume $X$ is connected.

Next if $\psi,\varphi\in W^{1,2}(X)\cap  L_\infty(X\!:\!\Ri)$ then by direct calculation
\[
\ci^{(h_N)}_\psi(\varphi)
=(\nabla(\varphi\psi),\nabla \psi)-2^{-1}(\nabla \varphi,\nabla(\psi^2))
=\int_Xdx\,\varphi(x) \, |(\nabla\psi)(x)|^2
\;\;\;.
\]
Therefore if $\psi\in W^{1,2}(X)_{\rm loc}\cap  L_\infty(X\!:\!\Ri)$ 
and  $|||\widehat \ci^{(h_N)}_\psi|||\leq 1$ one finds as in
 Example~\ref{xdfd318} that $\psi\in W^{1,\infty}(X\!:\!\Ri)$ and 
$\|\nabla\psi\|_\infty\leq 1$.
Hence $d^{(h_N)}(A\,;B)=\inf_{x\in A, y\in B}d(x\,;y)$ where
\[
d(x\,;y)=\sup\{|\psi(x)-\psi(y)|:W^{1,\infty}(X\!:\!\Ri)\, , \,\|\nabla\psi\|_\infty\leq 1\}
\;\;\;.
\]
But this is the geodesic distance, with  the usual Euclidean metric, from $x$ to $y$.
A similar conclusion follows for $d^{(h_D)}$ by replacing $W^{1,2}$ by $W^{1,2}_0$
in the argument. 
This replacement does not affect the  identification with the geodesic distance.
Therefore in both cases the set-theoretic distance between the sets is the geodesic distance 
in $X$ equipped with the Euclidean Riemannian metric.
\end{exam}

The definition of the distance in terms of the space $D(\ce)_{\rm loc}$ gives
good off-diagonal bounds but it is somewhat complicated. 
One could ask 
whether it has any simpler characterization in terms of $D(\ce)$.
One obvious choice is to set 
\[
d'_1(A\,;B)
= \sup \{ d_\psi(A\,;B) : \psi \in D_1(\ce)  \}
\;\;\;,
\]
for all measurable $A,B \subseteq X$ with  $\overline A, \overline B$ compact, where
\[
D_1(\ce)
= \{ \psi\in D(\ce) \cap L_\infty(X\!:\!\Ri) : |||\widehat \ci_\psi||| \leq 1\}
\;\;\;.
\]
Then one sets 
\[
d^{(\ce)}_1(A\,;B) = d_1(A\,;B)
= \inf \{ d'_1(A_0\,;B_0) : A_0\subseteq A, \; B_0\subseteq B \mbox{ and } \overline{A_0}, \overline{B_0} \mbox{ are compact} \}
\]
for all measurable $A,B \subseteq X$.
Note that $d_1(A\,;B) = d'_1(A\,;B)$ if $\overline A, \overline B$ are compact.
It follows that in general $d^{(\ce)}(A\,;B)\geq d^{(\ce)}_1(A\,;B)$ but one can have
a strict inequality.

\begin{exam} \label{xdfd512}
Let $\ce$ be the form associated with the second-order operator $-d^2/dx^2$ 
with Dirichlet boundary conditions on the interval $\langle0,1\rangle$. 
If $A=\langle0,a\rangle$ and $B=\langle b,1\rangle$ with $0<a<b<1$ then the boundary conditions
ensure that $d_\psi(A_\varepsilon\,;B_\varepsilon)\leq 2 \varepsilon$ for all $\psi\in D(\ce)\cap L_\infty(\langle0,1\rangle\!:\!\Ri)$
and small $\varepsilon > 0$, where $A_\varepsilon=\langle\varepsilon,a\rangle$ and $B_\varepsilon=\langle b,1-\varepsilon\rangle$.
Therefore $d_1(A\,;B)=0$.
But if $\psi(x) = 0 \vee (x-a) \wedge (b-a)$ then $\psi \in D_0(\ce)$
and $d(A\,;B) \geq d_\psi(A\,;B) = b-a$.
\end{exam}

The two distances are, however, often  equal.

\begin{prop}\label{pdfd513}
Let $\ce$ be a local Dirichlet form on $L_2(X)$.
Suppose that for all compact $K$ and $\varepsilon > 0$ there exists a $\chi \in D(\ce) \cap L_{\infty,c}(X\!:\!\Ri)$
such that $0 \leq \chi \leq 1$, $\chi|_K = 1$ and $|||\widehat \ci_\chi||| < \varepsilon$.
Then $d(A\,;B) = d_1(A\,;B)$ for all measurable $A,B \subseteq X$ with  $\overline A, \overline B$ compact.
\end{prop}
\proof\
We have to show that $d(A\,;B) \leq d_1(A\,;B)$.
Let $M \in [0,\infty\rangle$ and suppose that $M \leq d(A\,;B)$.
Let $\varepsilon \in \langle0,1]$.
By Lemma~\ref{ldfd438}.\ref{ldfd438-3} there exists a 
$\psi \in D_0(\ce)$ such that $d_\psi(A\,;B) \geq M - \varepsilon$
and $\|\psi\|_\infty \leq M$.
By assumption there exists a $\chi \in D(\ce) \cap L_{\infty,c}(X\!:\!\Ri)$
such that $0 \leq \chi \leq 1$, $\chi|_{\overline A \cup \overline B} = 1$ and $|||\widehat \ci_\chi||| < \varepsilon^2$.
Then $\chi \psi \in D(\ce) \cap L_\infty(X\!:\!\Ri)$ and it follows from Lemma~\ref{ldfd208}.\ref{ldfd208-2}
that 
\[
|||\widehat \ci_{\chi \psi}|||
\leq (1 + \delta) \, \|\chi\|_\infty^2 \, |||\widehat \ci_\psi|||
   + (1 + \delta^{-1}) \, \|\psi\|_\infty^2 \, |||\widehat \ci_\chi|||
\]
for all $\delta > 0$.
Choosing $\delta = \varepsilon$ gives
\[
|||\widehat \ci_{\chi \psi}|||
\leq (1 + \varepsilon) + (1 + \varepsilon^{-1}) M^2 \varepsilon^2
\leq 1 + \varepsilon(1 + 2 M^2)
\;\;\; .  \]
But $d_{\chi \psi}(A\,;B) = d_\psi(A\,;B) \geq M - \varepsilon$.
So $d_1(A\,;B) \geq (1 + \varepsilon(1 + 2 M^2))^{-1} (M - \varepsilon)$.
Hence $d_1(A\,;B) \geq M$ and the proposition follows.\hfill$\Box$

\section*{Acknowledgements}
This work  was 
supported by an Australian Research Council (ARC)
Discovery Grant DP 0451016.
Parts  of the collaboration took place during a visit of the
first, third and fourth named authors in  August 2004 to the
Australian National University and during a visit of the
second named author in  January 2005 to the
Eindhoven University of Technology.

We thank the referee for the reference to \cite{Schm}, which made it possible to 
remove several topological conditions on the space $X$ in the first version of this paper.


\begin{thebibliography}{ERSZ}

\bibitem[And]{And}
{\sc Andersson, L.-E.}, On the representation of Dirichlet forms.
\newblock {\em Ann.\ Inst.\ Fourier} {\bf 25} (1975),  11--25.

\bibitem[Aus]{Aus2}
{\sc Auscher, P.}, On necessary and sufficient conditions for $L^p$-estimates
  of Riesz transforms associated to elliptic operators on $\Ri^n$ and related
  estimates.
\newblock Research report, Preprint. Univ.\ de Paris-Sud, 2004.

\bibitem[BiM1]{BM1}
{\sc Biroli, M. {\rm and} Mosco, U.}, Formes de Dirichlet et estimations
  structurelles dans les milieux discontinus.
\newblock {\em C. R. Acad.\ Sci.\ Paris S\'er.\ I Math.} {\bf 313} (1991),
  593--598.

\bibitem[BiM2]{BM}
\leavevmode\vrule height 2pt depth -1.6pt width 23pt, A Saint-Venant type
  principle for Dirichlet forms on discontinuous media.
\newblock {\em Ann.\ Mat.\ Pura Appl.} {\bf 169} (1995).

\bibitem[Bra]{Bra}
{\sc Braides, A.}, {\em $\Gamma$-convergence for beginners}, vol.\ 22 of Oxford
  Lecture Series in Mathematics and its Applications.
\newblock Oxford University Press, Oxford, 2002.

\bibitem[BoH]{BH}
{\sc Bouleau, N. {\rm and} Hirsch, F.}, {\em Dirichlet forms and analysis on
  Wiener space}, vol.\ 14 of de Gruyter Studies in Mathematics.
\newblock Walter de Gruyter \& Co., Berlin, 1991.

\bibitem[CGT]{CGT}
{\sc Cheeger, J., Gromov, M. {\rm and} Taylor, M.}, Finite propagation speed,
  kernel estimates for functions of the Laplace operator, and the geometry of
  complete Riemannian manifolds.
\newblock {\em J. Differential Geom.} {\bf 17} (1982),  15--53.

\bibitem[{Dal}]{DalM}
{\sc {Dal Maso}, G.}, {\em An introduction to {$\Gamma$}-convergence}, vol.\ 8
  of Progress in Nonlinear Differential Equations and their Applications.
\newblock Birkh{\"a}user Boston Inc., Boston, MA, 1993.

\bibitem[Dav1]{Dav7}
{\sc Davies, E.B.}, Explicit constants for Gaussian upper bounds on heat
  kernels.
\newblock {\em Amer.\ J. Math.} {\bf 109} (1987),  319--333.

\bibitem[Dav2]{Dav12}
\leavevmode\vrule height 2pt depth -1.6pt width 23pt, Heat kernel bounds,
  conservation of probability and the Feller property.
\newblock {\em J. Anal.\ Math.} {\bf 58} (1992),  99--119.
\newblock Festschrift on the occasion of the 70th birthday of Shmuel Agmon.

\bibitem[EkT]{ET}
{\sc Ekeland, I. {\rm and} Temam, R.}, {\em Convex analysis and variational
  problems}.
\newblock North-Holland Publishing Co., Amsterdam, 1976.

\bibitem[ERSZ]{ERSZ1}
{\sc Elst, A.F.M. ter, Robinson, D.W., Sikora, A. {\rm and} Zhu, Y.},
  Second-order operators with degenerate coefficients.
\newblock Research Report CASA 04-32, Eindhoven University of Technology,
  Eindhoven, The Netherlands, 2004.

\bibitem[ERZ]{ERZ1}
{\sc Elst, A.F.M. ter, Robinson, D.W. {\rm and} Zhu, Y.}, Positivity and
  ellipticity.
\newblock {\em Proc.\ Amer.\ Math.\ Soc.} (2005).
\newblock To appear.

\bibitem[FOT]{FOT}
{\sc Fukushima, M., Oshima, Y. {\rm and} Takeda, M.}, {\em Dirichlet forms and
  symmetric Markov processes}, vol.\ 19 of de Gruyter Studies in Mathematics.
\newblock Walter de Gruyter \& Co., Berlin, 1994.

\bibitem[Gaf]{Gaf}
{\sc Gaffney, M.P.}, The conservation property of the heat equation on
  Riemannian manifolds.
\newblock {\em Comm.\ Pure Appl.\ Math.} {\bf 12} (1959),  1--11.

\bibitem[Gri]{Gri3}
{\sc Grigor'yan, A.}, Estimates of heat kernels on Riemannian manifolds.
\newblock In {\em Spectral theory and geometry $($Edinburgh, 1998$)$}, vol.\
  273 of London Math.\ Soc.\ Lecture Note Ser.,  140--225. Cambridge Univ.\
  Press, Cambridge, 1999.

\bibitem[HiR]{HiR}
{\sc Hino, M. {\rm and} Ram\'irez, J.A.}, Small-time Gaussian behavior of
  symmetric diffusion semigroups.
\newblock {\em Ann.\ Prob.} {\bf 31} (2003),  254--1295.

\bibitem[JeS1]{JSC}
{\sc Jerison, D.S. {\rm and} S{\'a}nchez-Calle, A.}, Estimates for the heat
  kernel for a sum of squares of vector fields.
\newblock {\em Ind.\ Univ.\ Math.\ J.} {\bf 35} (1986),  835--854.

\bibitem[JeS2]{JSC1}
\leavevmode\vrule height 2pt depth -1.6pt width 23pt, Subelliptic, second order
  differential operators.
\newblock In {\sc Berenstein, C.~A.}, ed., {\em Complex analysis III}, Lecture
  Notes in Mathematics 1277. Springer-Verlag, Berlin etc., 1987,  46--77.

\bibitem[Jos]{Jos}
{\sc Jost, J.}, Nonlinear Dirichlet forms.
\newblock In {\em New directions in Dirichlet forms}, vol.\ 8 of AMS/IP Stud.\
  Adv.\ Math.,  1--47. Amer.\ Math.\ Soc., Providence, RI, 1998.

\bibitem[Kat]{Kat1}
{\sc Kato, T.}, {\em Perturbation theory for linear operators}.
\newblock Second edition, Grundlehren der mathematischen Wissenschaften 132.
  Springer-Verlag, Berlin etc., 1984.

\bibitem[MaR]{MR}
{\sc Ma, Z.M. {\rm and} R{\"o}ckner, M.}, {\em Introduction to the theory of
  (non symmetric) Dirichlet Forms}.
\newblock Universitext. Springer-Verlag, Berlin etc., 1992.

\bibitem[Mos]{Mosco}
{\sc Mosco, U.}, Composite media and asymptotic Dirichlet forms.
\newblock {\em J. Funct.\ Anal.} {\bf 123} (1994),  368--421.

\bibitem[ReS]{RS4}
{\sc Reed, M. {\rm and} Simon, B.}, {\em Methods of modern mathematical physics
  IV. Analysis of operators}.
\newblock Academic Press, New York etc., 1978.

\bibitem[Rot]{Roth}
{\sc Roth, J.-P.}, Formule de repr\'esentation et troncature des formes de
  Dirichlet sur $\Ri^m$.
\newblock In {\em S\'eminaire de Th\'eorie du Potentiel de Paris, No. 2},
  Lect.\ Notes in Math. 563,  260--274. Springer Verlag, Berlin, 1976.

\bibitem[Sch]{Schm}
{\sc Schmuland, B.}, On the local property for positivity preserving coercive
  forms.
\newblock In {\sc Ma, Z.M. {\rm and} R{\"o}ckner, M.}, eds., {\em Dirichlet
  forms and stochastic processes}. Walter de Gruyter \& Co., Berlin, 1995,
  345--354.
\newblock Papers from the International Conference held in Beijing, October
  25--31, 1993, and the School on Dirichlet Forms, held in Beijing, October
  18--24, 1993.

\bibitem[Sik]{Sik3}
{\sc Sikora, A.}, Riesz transform, Gaussian bounds and the method of wave
  equation.
\newblock {\em Math.\ Z.} {\bf 247} (2004),  643--662.

\bibitem[Sim1]{bSim6}
{\sc Simon, B.}, Ergodic semigroups of positivity preserving self-adjoint
  operators.
\newblock {\em J.\ Funct.\ Anal.} {\bf 12} (1973),  335--339.

\bibitem[Sim2]{bSim4}
\leavevmode\vrule height 2pt depth -1.6pt width 23pt, Lower semicontinuity of
  positive quadratic forms.
\newblock {\em Proc.\ Roy.\ Soc.\ Edinburgh Sect.\ A} {\bf 79} (1977),
  267--273.

\bibitem[Sim3]{bSim5}
\leavevmode\vrule height 2pt depth -1.6pt width 23pt, A canonical decomposition
  for quadratic forms with applications to monotone convergence theorems.
\newblock {\em J.\ Funct.\ Anal.} {\bf 28} (1978),  377--385.

\bibitem[Stu1]{Stu4}
{\sc Sturm, K.-T.}, Analysis on local Dirichlet spaces. II. Upper Gaussian
  estimates for the fundamental solutions of parabolic equations.
\newblock {\em Osaka J. Math.} {\bf 32} (1995),  275--312.

\bibitem[Stu2]{Stu2}
\leavevmode\vrule height 2pt depth -1.6pt width 23pt, The geometric aspect of
  Dirichlet forms.
\newblock In {\em New directions in Dirichlet forms}, vol.\ 8 of AMS/IP Stud.\
  Adv.\ Math.,  233--277. Amer.\ Math.\ Soc., Providence, RI, 1998.

\end{thebibliography}
\end{document}